\documentclass[reqno,pdfmx,dvipdfmx]{amsart}

\PassOptionsToPackage{dvipdfmx}{graphicx}
\usepackage{amsmath,amsthm,amssymb,graphicx,cite,bm,mathrsfs,latexsym,hyperref}

\theoremstyle{plain}
\newtheorem{thm}{Theorem}[section]

\newtheorem{rmk}[thm]{Remark}

\def\D{\mathrm{D}}

\def\I{\mathrm{I}}
\def\R{\mathrm{R}}

\def\T{\mathrm{T}}

\def\d{\mathrm{d}}

\def\i{\mathrm{i}}

\def\s{\mathrm{s}}
\def\u{\mathrm{u}}

\def\Cset{\mathbb{C}}
\def\Nset{\mathbb{N}}
\def\Rset{\mathbb{R}}

\def\Zset{\mathbb{Z}}

\def\epsilon{\varepsilon}

\def\e{\mathrm{e}}

\DeclareMathOperator{\sech}{sech}
\DeclareMathOperator{\gKer}{gKer}
\DeclareMathOperator{\Ker}{Ker}

\DeclareMathOperator{\sgn}{sgn}

\renewcommand{\Re}{\mathrm{Re}\,}
\renewcommand{\Im}{\mathrm{Im}\,}

\newcommand\defeq{\mathrel{\rlap{\raisebox{0.3ex}{$\cdot$}}\raisebox{-0.3ex}{$\cdot$}}=}

\newcommand*\pFqskip{8mu}
\catcode`,\active
\newcommand*\pFq{\begingroup
        \catcode`\,\active
        \def ,{\mskip\pFqskip\relax}%
        \dopFq
}
\catcode`\,12
\def\dopFq#1#2#3#4#5{%
        {}_{#1}F_{#2}\biggl(\genfrac..{0pt}{}{#3}{#4};#5\biggr)%
        \endgroup
}

\makeatletter
 \@addtoreset{equation}{section}
\makeatother

\begin{document}


\title[]{
Numerical computations for bifurcations and spectral stability
of solitary waves in coupled nonlinear Schr\"odinger equations}
\thanks{This work was partially supported by JSPS KAKENHI Grant Number JP17H02859.}

\author{Kazuyuki Yagasaki}

\author{Shotaro Yamazoe}

\address{Department of Applied Mathematics and Physics, Graduate School of Informatics,
Kyoto University, Yoshida-Honmachi, Sakyo-ku, Kyoto 606-8501, JAPAN}
\email{yagasaki@amp.i.kyoto-u.ac.jp (K. Yagasaki)}
\email{yamazoe@amp.i.kyoto-u.ac.jp (S. Yamazoe)}

\date{\today}

\begin{abstract}
We numerically study solitary waves in the coupled nonlinear Schr\"odinger equations.
We detect pitchfork bifurcations of the \emph{fundamental solitary wave}
 and compute eigenvalues and eigenfunctions of the corresponding eigenvalue problems
 to determine the spectral stability of solitary waves born at the pitchfork bifurcations.
Our numerical results demonstrate the theoretical ones which the authors obtained recently.
We also compute generalized eigenfunctions associated with the zero eigenvalue
 for the bifurcated solitary wave exhibiting a saddle-node bifurcation,
 and show that it does not change its stability type
 at the saddle-node bifurcation point.
\end{abstract}
\maketitle


\section{Introduction}

We consider the coupled nonlinear Schr\"odinger (CNLS) equations of the form
\begin{align}
	\left.
		\begin{array}{l}
			\i\partial_t u=-\partial_x^2u-(|u|^2+\beta_1|v|^2)u,\\
			\i\partial_t v=-\partial_x^2v-(\beta_1|u|^2+\beta_2|v|^2)v,
		\end{array}
	\right.
	\quad(t,x)\in\Rset\times\Rset,
	\label{eqn:CNLS}
\end{align}
where $(u,v) = (u(t,x),v(t,x))$ are complex-valued unknown functions of $(t,x)\in \Rset\times\Rset$ and $\beta_1,\beta_2\in\Rset$ are parameters.
Here we are interested in the solitary wave solutions to \eqref{eqn:CNLS} of the form
\begin{align}
	\begin{aligned}
		&u(t,x)=\e^{\i(\omega t+cx-c^2t+\theta)}U(x-2ct-x_0),\\
		&v(t,x)=\e^{\i(st+cx-c^2t+\phi)}V(x-2ct-x_0),
	\end{aligned}
	\label{eqn:soliton}
\end{align}
where $\omega,s>0$, and $c,x_0,\theta,\phi\in\Rset$ are constants, such that the real-valued functions $(U,V)=(U(x),V(x))$ satisfy $U(x),V(x)\to 0$ as $x\to\pm\infty$.
Henceforth, without loss of generality, we take $c,x_0,\theta,\phi=0$ since \eqref{eqn:CNLS} is invariant under the Galilean transformations
\begin{align*}
	(u(t,x),v(t,x))\mapsto\e^{\i(cx-c^2t)}(u(t,x-2ct),v(t,x-2ct)),\quad c\in\Rset,
\end{align*}
the spatial translations
\begin{align*}
	(u(t,x),v(t,x))\mapsto(u(t,x-x_0),v(t,x-x_0)),\quad x_0\in\Rset,
\end{align*}
and the gauge transformations
\begin{align*}
	(u(t,x),v(t,x))\mapsto(\e^{\i\theta}u(t,x),\e^{\i\phi}v(t,x)),\quad\theta,\phi\in\Rset.
\end{align*}
So $(U,V)=(U(x),V(x))$ solves
\begin{align}
	\begin{aligned}
		&-U''+\omega U-(U^2+\beta_1V^2)U=0,\\
		&-V''+sV-(\beta_1U^2+\beta_2V^2)V=0,
	\end{aligned}
	\label{eqn:SE}
\end{align}
where the prime represents the differentiation with respect to $x$.
In particular, \eqref{eqn:SE} allows homoclinic solutions
 of which one component is identically zero, e.g.,
\begin{equation}
(U,V)=(U_0(x),0),\quad
U_0(x)\defeq\sqrt{2\omega}\sech(\sqrt\omega x).
\label{eqn:homo}
\end{equation}
We refer to the solitary waves corresponding
 to such homoclinic solutions in \eqref{eqn:CNLS}
 as the \textit{fundamental solitary waves}.
Bl\'azquez-Sanz and Yagasaki \cite{BY12} showed that
 the homoclinic solution $(U_0(x),0)$ exhibits infinitely many pitchfork bifurcations
 in \eqref{eqn:SE} when $\beta_1$ is increased from zero for $\beta_2$ fixed.
This result means that the fundamental solitary wave
\begin{align}
	(u(t,x),v(t,x))=(\e^{\i\omega t}U_0(x),0)
	\label{eqn:fun_sw}
\end{align}
also exhibits infinitely many ones in \eqref{eqn:CNLS}.
Here the terminology ``pitchfork bifurcation'' is used with caution:
 a pair of homoclinic solutions to \eqref{eqn:SE},
 which correspond to the same family of solitary waves of the form \eqref{eqn:soliton}
 in \eqref{eqn:CNLS}, are born at the bifurcation point.

Bifurcations and stability of solitary waves in nonlinear wave equations
 have been widely investigated \cite{KaP13,P11,Ya10}.
For CNLS equations with the cubic nonlinearity, internal oscillations and radiation dumping of the single-hump vector solitons were studied by Yang \cite{Ya97b} and Pelinovsky and Yang \cite{PYa00}.
For general nonlinearity cases, Yang \cite{Ya97a} classified possible bifurcations of solitary waves, and Pelinovsky and Yang \cite{PYa05} determined the stability of solitary waves under some generic nondegenerate conditions.
Jackson \cite{J11} also studied the stability of solitary waves
 from a geometric point of view.
Recently, the authors \cite{YY20b}
 used the approach of \cite{BY12} and developed some techniques
 to detect pitchfork bifurcations of the fundamental solitary wave
 and spectral stability of the fundamental and bifurcated solitary waves
 in general CNLS equations containing \eqref{eqn:CNLS}.
A perturbation expansion of eigenvalues of the linearized operator
 was directly calculated under some nondegenerate conditions
 which are easy to verify compared to assumptions made in \cite{PYa05}.
In particular, for \eqref{eqn:CNLS},
 it was shown in \cite{YY20b} that
 the solitary waves born at the first bifurcation are stable
 but the solitary waves born at  the other bifurcations are unstable
 while the fundamental one continues to be stable.

In this paper, we numerically detect the pitchfork bifurcations
 of the fundamental solitary wave \eqref{eqn:fun_sw}
 and compute eigenvalues and eigenfunctions of the corresponding eigenvalue problems
 to determine the spectral stability of solitary waves born at the pitchfork bifurcations
 in \eqref{eqn:CNLS}.
In particular, the numerical results, some of which were also provided in \cite{YY20b},
 demonstrate the theoretical ones obtained in \cite{YY20b}.
Moreover, one of the bifurcated solitary waves is observed
 to exhibit a saddle-node bifurcation.
We compute generalized eigenfunctions associated with the zero eigenvalue
 to show that the solitary wave does not change its stability type
 at the saddle-node bifurcation point.
Such saddle-node bifurcations with no stability switching
 were proven to occur in general single nonlinear Schr\"odinger (NLS) equations
 with external potentials by Yang \cite{Ya12b} earlier.
 To the authors' knowledge, such a phenomenon has not been reported for CNLS equations before.
The computer tool \texttt{AUTO} \cite{D12} was used
 for carrying out necessary computations,
 as in similar numerical work of \cite{YY20a}
 for the single NLS equation with an external potential.

This paper is organized as follows:
In Section~2 we briefly review the theoretical results of \cite{YY20b}
 on bifurcations and the stability of solitary waves in \eqref{eqn:CNLS}.
We give numerical computations for eigenvalues and eigenfunctions in Section~3
 and for generalized eigenfunctions associated with the zero eigenvalues in Section~4.
Our numerical approaches, of which a general framework was given in \cite{YY20a} for eigenvalues and eigenfunctions,
 are briefly described there before the results are provided.


\section{Theoretical Results}

In this section we briefly review the theoretical results of \cite{YY20b}
 on bifurcations of the fundamental solitary wave \eqref{eqn:fun_sw} 
 and the stability of the fundamental and bifurcated solitary waves
 in the CNLS equations \eqref{eqn:CNLS}.

We begin with the bifurcation result.
The variational equation (VE) of \eqref{eqn:SE}
 around the homoclinic solution $(U,V)=(U_0(x),0)$ is given by
\begin{equation}
	-\delta U''+\omega\delta U-3U_0(x)^2\delta U=0,\quad
	-\delta V''+s\delta V-\beta_1U_0(x)^2\delta V=0.
	\label{eqn:VE}
\end{equation}
We easily see that $(\delta U,\delta V)=(U_0'(x),0)$
 is a bounded solution to \eqref{eqn:VE}.
We also show that \eqref{eqn:VE} has another bounded solution
 $(\delta U,\delta V)=(0,V_1^{(\ell)}(x))$, which is linearly independent of $(U_0'(x),0)$,
 if and only if
\begin{align}
	\beta_1=\beta_1^{(\ell)}\defeq\frac{(\sqrt{s/\omega}+\ell)(\sqrt{s/\omega}+\ell+1)}{2},\quad
	\ell\in\Zset_{\ge 0}\defeq\{n\in\Zset\mid n\ge 0\},
	\label{eqn:beta1_ell}
\end{align}
where
\begin{align*}
	V_1^{(\ell)}(x)=
	\sech^{\sqrt{s/\omega}}(\sqrt\omega x)\,\pFq{2}{1}{-\ell',\sqrt{s/\omega}+\ell'+1/2}{\sqrt{s/\omega}+1}{\sech^2(\sqrt\omega x)}
\end{align*}
for $\ell=2\ell'$, $\ell'\in\Zset_{\ge 0}$, and
\begin{align*}
	V_1^{(\ell)}(x)=
	\sech^{\sqrt{s/\omega}}(\sqrt\omega x)\tanh(\sqrt\omega x)\,\pFq{2}{1}{-\ell',\sqrt{s/\omega}+\ell'+3/2}{\sqrt{s/\omega}+1}{\sech^2(\sqrt\omega x)}
\end{align*}
for $\ell=2\ell'+1$, $\ell'\in\Zset_{\ge 0}$.
Here ${}_2F_1$ denotes the hypergeometric function
\begin{align*}
	\pFq{2}{1}{a,b}{c}{z}=\sum_{j=0}^\infty\frac{(a)_j(b)_j}{j!\,(c)_j}z^j,
\end{align*}
where $a,b,c$ are constants, $(x)_j\defeq\Gamma(x+j)/\Gamma(x)$,
 and $\Gamma(x)$ is the Gamma function.
See Section~5 of \cite{BY12} or Section~6 of \cite{YY20b}.
The case of $\omega=1$ was considered there
 by replacing $\omega t$ and $s/\omega$ with $t$ and $s$, respectively,
 without loss of generality.

Define the integral,
\begin{align*}
&
\bar a_2=-2\int_{-\infty}^\infty V_1^{(\ell)}(x)^2\sech^2x\,\d x<0,\\
&
\bar b_2=8\bigl(\beta_1^{(\ell)}\bigr)^2\int_{-\infty}^\infty\phi_{11}(x)V_1^{(\ell)}(x)^2\sech x\,\biggl(\int_x^\infty\phi_{12}(y)V_1^{(\ell)}(y)^2\sech y\,\d y\biggr)\d x\notag\\
	&\qquad-\beta_2\int_{-\infty}^\infty V_1^{(\ell)}(x)^4\,\d x,
\end{align*}
where
\begin{align*}
	\phi_{11}(x)=\frac{1}{2}\sech x\,(3-\cosh^2x-3x\tanh x),\quad
	\phi_{12}(x)=\sech x\tanh x.
\end{align*}
We note that $\phi_{11}(x)$ and $\phi_{12}(x)$ are, respectively,
 the $(1,1)$- and $(1,2)$-elements of a fundamental matrix $\Phi(x)$
 of the linear system
\[
\frac{\d}{\d x}
\begin{pmatrix}
\delta U\\
\delta U'
\end{pmatrix}
=\begin{pmatrix}
0 & 1\\
\omega-3U_0(x)^2 & 0
\end{pmatrix}
\begin{pmatrix}
\delta U\\
\delta U'
\end{pmatrix},
\]
as which the first equation of the VE \eqref{eqn:VE} is rewritten in a first-order system,
 such that $\Phi(0)=I_2$,
 where $I_n$ denotes the $n\times n$ identity matrix for $n>1$.
The following result was proven
 on bifurcations of the fundamental solitary wave \eqref{eqn:fun_sw}
 in Theorem~7.1 of \cite{YY20b} (see also Theorem~5.3 (ii) of \cite{BY12}).

\begin{thm}\label{thm:bif}
	For $\ell\in\Zset_{\ge 0}$,
	a pitchfork bifurcation of the fundamental solitary wave \eqref{eqn:fun_sw}
	occurs at $\beta_1=\beta_1^{(\ell)}$ if $\bar{b}_2\neq 0$.
	In addition, it is supercritical or subcritical,
	depending on whether $\bar{b}_2>0$ or $<0$.
	Moreover, the bifurcated solitary waves are expressed as
	\begin{align}
		(u(t,x),v(t,x))=(\e^{\i\omega t}U_\epsilon(x),\e^{\i st}V_\epsilon(x))
		\label{eqn:bif_sw}
	\end{align}
	with
	\begin{align}
		U_\epsilon(x)=U_0(x)+O(\epsilon^2),\quad
		V_\epsilon(x)=\epsilon V_1^{(\ell)}(x)+O(\epsilon^3),
	\label{eqn:UVe}
	\end{align}
	where $\epsilon>0$ is a small parameter such that $\beta_1=\beta_1^{(\ell)}+O(\epsilon^2)$.
\end{thm}

A more precise expression of the bifurcated solitary waves than \eqref{eqn:UVe}
 was given in Theorem~7.1 of \cite{YY20b} (see also Theorem~2.2 of \cite{YY20b}).
Tractable expressions of the integrals $\bar{a}_2$ and $\bar{b}_2$
 for computation were also obtained in Proposition~7.4 of \cite{YY20b}
 (see also Appendix~B of \cite{YY20b}
  for closed-form ones of $\bar{b}_2$ when $\ell\le4$).

We turn to the stability result.
The linearized operator of \eqref{eqn:CNLS}
 around the solitary wave \eqref{eqn:soliton} with $c,x_0,\theta,\phi=0$
 is given by $J\mathcal L$ with
\begin{align}
	J\defeq\begin{pmatrix}O_2&I_2\\-I_2&O_2\end{pmatrix},\quad
	\mathcal L\defeq\begin{pmatrix}\mathcal L_+&O_2\\O_2&\mathcal L_-\end{pmatrix},
\label{eqn:JL}
\end{align}
where $O_n$ is the $n\times n$ zero matrix for $n\in\Nset$ and $n>1$, and
\begin{align}
	\begin{aligned}
		&\mathcal L_+\defeq\begin{pmatrix}
			-\partial_x^2+\omega-(3U^2+\beta_1V^2)&-2\beta_1UV\\
			-2\beta_1UV&-\partial_x^2+s-(\beta_1U^2+3\beta_2V^2)
		\end{pmatrix},\\
		&\mathcal L_-\defeq\begin{pmatrix}
			-\partial_x^2+\omega-(U^2+\beta_1V^2)&0\\
			0&-\partial_x^2+s-(\beta_1U^2+\beta_2V^2)
		\end{pmatrix}.
	\end{aligned}
	\label{eqn:def_Lpm}
\end{align}
See Section~4.2 of \cite{YY20b}
 for the derivation of the expression of $J\mathcal L$
 in more general CNLS equations containing \eqref{eqn:CNLS}.
To discuss the spectral stability of the solitary wave,
 we consider the associated eigenvalue problem
\begin{align}
	J\mathcal L\psi=\lambda\psi,\quad
	\psi\in L^2(\Rset)^4.
	\label{eqn:eig_eq}
\end{align}
We easily obtain the following properties of the spectrum $\sigma(J\mathcal L)$
 (see Section~4.1 of \cite{YY20b} for the details):
\begin{enumerate}
\setlength{\leftskip}{-2em}
\item[(i)]
If $\lambda\in\sigma(J\mathcal L)$,
 then $-\lambda,\pm\overline\lambda\in\sigma(J\mathcal L)$,
 where the overline represents the complex conjugate.
Actually, if $\psi=(\psi_1,\ldots,\psi_4)^\T$
 is an eigenfunction of \eqref{eqn:eig_eq} for the eigenvalue $\lambda$, then
\begin{align*}
	J\mathcal L\overline\psi=\overline\lambda\overline\psi,\quad
	J\mathcal L\begin{pmatrix}\psi_1\\\psi_2\\-\psi_3\\-\psi_4\end{pmatrix}
	=-\lambda\begin{pmatrix}\psi_1\\\psi_2\\-\psi_3\\-\psi_4\end{pmatrix}.
\end{align*}
\item[(ii)]
The essential spectrum is given by
\begin{align}
	\sigma_{\mathrm{ess}}(J\mathcal L)=\i(-\infty,-\min\{\omega,s\}]\cup\i[\min\{\omega,s\},\infty).\label{eqn:ess_spec}
\end{align}
\item[(iii)]
$\Ker J\mathcal L=\Ker\mathcal L_+\oplus\Ker\mathcal L_-$ contains
\begin{align}
	\begin{aligned}
		&\varphi_1(x)=(U'(x),V'(x),0,0)^\T,\\
		&\varphi_2(x)=(0,0,U(x),0)^\T,\quad
		\varphi_3(x)=(0,0,0,V(x))^\T.
	\end{aligned}
	\label{eqn:kernel}
\end{align}
Moreover,
\begin{align}
	\begin{aligned}
		&\chi_1(x)=(0,0,-xU(x)/2,-xV(x)/2)^\T,\\
		&\chi_2(x)=(\partial_\omega U(x),\partial_\omega V(x),0,0)^\T,\quad
		\chi_3(x)=(\partial_s U(x),\partial_s V(x),0,0)^\T
	\end{aligned}
	\label{eqn:def_chi}
\end{align}
satisfy $J\mathcal L\chi_j=\varphi_j$, $j=1,2,3$, whenever they exist.
\end{enumerate}
Let $\beta_1>0$ and let $\kappa=(-1+\sqrt{1+8\beta_1})/2$.
For $(U,V)=(U_0,0)$ (i.e., the fundamental solitary wave \eqref{eqn:fun_sw}), 
 $J\mathcal L$ has the eigenvalues
\begin{equation}
\lambda=\pm\i(s-\omega(\kappa-k)^2),\quad
k\in\{0,1,\ldots,\lfloor\kappa\rfloor\}\setminus\{\kappa\},
\label{eqn:ev}
\end{equation}
and the associated eigenfunctions $\psi=(0,\Psi,0,\pm\i\Psi)^\T$ with
\begin{align}
	\Psi(x)=\sech^{\kappa-k}(\sqrt\omega x)\,\pFq{2}{1}{-k',\kappa-k+k'+1/2}{\kappa-k+1}{\sech^2(\sqrt\omega x)}
\label{eqn:ef1}
\end{align}
if $k\ne\kappa$ is even ($k=2k'$, $k'\in\Zset_{\ge 0}$), and with
\begin{align}
	\Psi(x)=\sech^{\kappa-k}(\sqrt\omega x)\tanh(\sqrt\omega x)\,\pFq{2}{1}{-k',\kappa-k+k'+3/2}{\kappa-k+1}{\sech^2(\sqrt\omega x)}
\label{eqn:ef2}
\end{align}
if $k\ne\kappa$ is odd ($k=2k'+1$, $k'\in\Zset_{\ge 0}$),
 where the upper or lower signs are taken simultaneously.
See Remark~7.6 of \cite{YY20b}.

By analyzing the eigenvalue problem \eqref{eqn:eig_eq}
 based on the Evans function technique \cite{AGJ90,LP00},
 the following result was proven in Theorem~7.9 of \cite{YY20b}.

\begin{thm}
\label{thm:ls}
	The solitary wave \eqref{eqn:bif_sw} born at $\beta_1=\beta_1^{(0)}$ near there
	as well as the fundamental solitary wave \eqref{eqn:fun_sw} for $\beta_1\in\Rset$
	is spectrally and orbitally stable.
	The solitary wave \eqref{eqn:bif_sw} born at $\beta_1=\beta_1^{(\ell)}$ with $\ell\ge1$
	is spectrally unstable if
	\begin{align}
		\sqrt{s/\omega}\notin\left\{-\frac{\ell^2+k^2}{2(\ell+k)}\;\middle|\;k\in\{-\lfloor(\sqrt 2+1)\ell\rfloor,\ldots,-\ell-1\}\right\}.
		\label{eqn:thmls}
	\end{align}
\end{thm}

\begin{rmk}
{\setlength{\leftmargini}{1.8em}
\begin{enumerate}
\item[(i)]
	If condition~\eqref{eqn:thmls} does not hold,
	 then some purely imaginary eigenvalues of $J\mathcal L$
	 around the fundamental solitary wave are of multiplicity two,
	 so that further tremendous treatments are required to determine their stability.
\item[(ii)]
	Pelinovsky and Yang \cite{PYa05} obtained a similar result under some generic conditions
	 which are difficult to actually check for \eqref{eqn:CNLS}.
\item[(iii)]
	The mechanism of instability for $\ell\ge1$ is stated as follows.
	The eigenvalues \eqref{eqn:ev} are embedded in the essential spectrum \eqref{eqn:ess_spec}.
	Moreover, they have a negative Krein signature if $0\le k<\ell$.
	Since eigenvalues with a negative Krein signature are structually unstable \cite{G90}, they split to a pair of eigenvalues with positive and negative real parts under perturbations generically.
\end{enumerate}
}
\end{rmk}



\section{Computations of eigenvalues and eigenfunctions}

In this section,
 we give some numerical computation results for eigenvalues and eigenfunctions
 of the eigenvalue problem \eqref{eqn:eig_eq}
 along with homoclinic solutions to \eqref{eqn:SE},
 and demonstrate the theoretical results stated in Section~2
 on bifurcations of the fundamental solitary wave \eqref{eqn:fun_sw}
 and the stability of bifurcated solitary waves in the CNLS equations \eqref{eqn:CNLS}
 by the numerical ones.


\subsection{Numerical Approach}
We first briefly describe our numerical approach,
 which was provided in a general setting in Section~2 of \cite{YY20a}.

We begin with computation of homoclinic solutions to \eqref{eqn:SE}.
We slightly modify \eqref{eqn:SE}, 
 rewrite it in a first-order system as
\begin{align}
	z'=f(z;\beta_1),\quad z=(z_1,z_2,z_3,z_4)^\T\in\Rset^4,
	\label{eqn:dyn_sys}
\end{align}
and numerically compute a homoclinic solution to \eqref{eqn:dyn_sys} satisfying
\begin{align}
	\lim_{x\to\pm\infty}z(x)=0,
	\label{eqn:z_decay}
\end{align}
where
\begin{align*}
	f(z;\beta_1)
	\defeq\begin{pmatrix}
		z_3\\
		z_4\\
		\omega z_1-(z_1^2+\beta_1z_2^2)z_1+d_1z_3\\
		sz_2-(\beta_1z_1^2+\beta_2z_2^2)z_2+d_1z_4
	\end{pmatrix}
\end{align*}
with a dummy parameter $d_1$.
Note that 
 \eqref{eqn:dyn_sys} is equivalent to \eqref{eqn:SE} if $d_1=0$.
We perform continuation of homoclinic orbits with two parameters since their existence is of codimension one \cite{CKS96}.
Moreover, a homoclinic solution persists in \eqref{eqn:dyn_sys}
 when one of the other parameters changes only if $d_1=0$
 (see Lemma~2.13 and Section~5 of \cite{BY12}).

Let $E^\s$ and $E^\u$ be, respectively,
 the two-dimensional stable and unstable subspaces
 of the linearized system at the origin for \eqref{eqn:dyn_sys},
\begin{align}
	\delta z'=\D_zf(0;\beta_1)\delta z.
\label{eqn:dyn_sys0}
\end{align}
We approximate the homoclinic solution $z(x)$ to \eqref{eqn:dyn_sys}
 satisfying \eqref{eqn:z_decay},
 so that it starts on $E^\u$ near the origin at $x_-$
 and arrives on $E^\s$ near the origin at $x_+$,
 where $x_-<0<x_+$ and $|x_\pm|\gg1$.
So we look for a solution to \eqref{eqn:dyn_sys} satisfying
\begin{align}
	L^\s z(x_-)=0,\quad
	L^\u z(x_+)=0,
	\label{eqn:bc1}
\end{align}
where
\begin{align*}
	L^\s\defeq\begin{pmatrix}
		-d_1/2-\sqrt{\omega+(d_1/2)^2}&0&1&0\\
		0&-d_1/2-\sqrt{s+(d_1/2)^2}&0&1
	\end{pmatrix}
\end{align*}
and
\begin{align*}
	L^\u\defeq\begin{pmatrix}
		-d_1/2+\sqrt{\omega+(d_1/2)^2}&0&1&0\\
		0&-d_1/2+\sqrt{s+(d_1/2)^2}&0&1
	\end{pmatrix}
\end{align*}
are $2\times4$ real matrices consisting of row eigenvectors
 for $\D_zf(0;\beta_1)$
 such that the associated eigenvalues are negative and positive, respectively.
The distances $|z(x_\pm)|$ should be kept small in the computation.
To eliminate the multiplicity of solutions due to the translational symmetry of \eqref{eqn:SE}, we also add the integral condition
\begin{align}
	\sum_{j=1}^2\int_{x_-}^{x_+}\bigl(z_j(x)-z_j^*(x)\bigr)z_{j+2}^*(x)\,\d x=0,
	\label{eqn:ic1}
\end{align}
where $z^*=(z^*_1,\ldots,z^*_4)^\T$ represents a previously computed solution
 along a continuation branch.

We turn to the eigenvalue problem \eqref{eqn:eig_eq} and rewrite it as
\begin{align}
	\zeta'=A(x;\beta_1,\lambda)\zeta,\quad
	\zeta\in\Cset^8,\ \lambda\in\Cset,
	\label{eqn:eig_eq_2}
\end{align}
with
\begin{align}
	\lim_{x\to\pm\infty}\zeta(x)=0,
	\label{eqn:zeta_decay}
\end{align}
where
\begin{align*}
A(x;\beta_1,\lambda)\defeq
\left(\begin{array}{c|c}
O_4&I_4\\\hline
\begin{array}{cc}
A_1(x;\beta_1,\lambda) & \lambda I_2\\
-\lambda I_2 & A_2(x;\beta_1,\lambda)\\
\end{array}
&O_4
\end{array}\right)
\end{align*}
with
\begin{align*}
&
A_1(x;\beta_1,\lambda)
\defeq\left(\begin{array}{cc}
\omega-(3U^2+\beta_1V^2)&-2\beta_1UV\\[1ex]
-2\beta_1UV&s-(\beta_1U^2+3\beta_2V^2)
\end{array}\right),\\
&
A_2(x;\beta_1,\lambda)
\defeq\left(\begin{array}{cc}
\omega-(U^2+\beta_1V^2)&0\\[1ex]
0&s-(\beta_1U^2+\beta_2V^2)
\end{array}\right).
\end{align*}
Let $\lambda=\lambda_\R+\i\lambda_\I$ with $\lambda_\R,\lambda_\I\in\Rset$ and let
\begin{align*}
	A(x;\beta_1,\lambda)=A_\R(x;\beta_1,\lambda_\R,\lambda_\I)+\i A_\I(x;\beta_1,\lambda_\R,\lambda_\I),
\end{align*}
where $A_\R(x;\beta_1,\lambda_\R,\lambda_\I)$ and $A_\I(x;\beta_1,\lambda_\R,\lambda_\I)$ are $8\times8$ real matrices.
Letting $\zeta=\zeta_\R+\i\zeta_\I$ with $\zeta_\R,\zeta_\I\in\Rset^8$,
 we rewrite \eqref{eqn:eig_eq_2} and \eqref{eqn:zeta_decay} as
\begin{equation}
	\begin{split}
		&\zeta_\R'
		=A_\R(x;\beta_1,\lambda_\R,\lambda_\I)\zeta_\R-A_\I(x;\beta_1,\lambda_\R,\lambda_\I)\zeta_\I,\\
		&\zeta_\I'
		=A_\I(x;\beta_1,\lambda_\R,\lambda_\I)\zeta_\R+A_\R(x;\beta_1,\lambda_\R,\lambda_\I)\zeta_\I
	\end{split}
	\label{eqn:eig_eq_2'}
\end{equation}
and
\begin{align}
	\lim_{x\to\pm\infty}\zeta_\R(x)=\lim_{x\to\pm\infty}\zeta_\I(x)=0,
	\label{eqn:bc2'}
\end{align}
respectively.

Let
\begin{align*}
	A_\infty(\beta_1,\lambda)
	\defeq\lim_{x\to\pm\infty}A(x;\beta_1,\lambda)
	=\left(
		\begin{array}{c|c}
			O_4&I_4\\\hline
			\begin{array}{cccc}
				\omega&0&\lambda&0\\
				0&s&0&\lambda\\
				-\lambda&0&\omega&0\\
				0&-\lambda&0&s
			\end{array}&O_4
		\end{array}
	\right)
\end{align*}
and let $\tilde{E}^\s$ and $\tilde{E}^\u$ be, respectively,
 the four-dimensional stable and unstable subspaces of the autonomous linear system
\begin{equation*}
	\begin{split}
		&\zeta_\R'
		=A_{\R\infty}(\beta_1,\lambda_\R,\lambda_\I)\zeta_\R-A_{\I\infty}(\beta_1,\lambda_\R,\lambda_\I)\zeta_\I,\\
		&\zeta_\I'
		=A_{\I\infty}(\beta_1,\lambda_\R,\lambda_\I)\zeta_\R+A_{\R\infty}(\beta_1,\lambda_\R,\lambda_\I)\zeta_\I,
	\end{split}
\end{equation*}
where $A_{\R\infty}(\beta_1,\lambda_\R,\lambda_\I)$ and $A_{\I\infty}(\beta_1,\lambda_\R,\lambda_\I)$ are $8\times8$ real matrices such that
\begin{align*}
	A_\infty(\beta_1,\lambda)=A_{\R\infty}(\beta_1,\lambda_\R,\lambda_\I)+\i A_{\I\infty}(\beta_1,\lambda_\R,\lambda_\I).
\end{align*}
Like the homoclinic solution to \eqref{eqn:dyn_sys},
 we approximate the solution $(\zeta_\R(x),\zeta_\I(x))$ to \eqref{eqn:eig_eq_2'}
 satisfying \eqref{eqn:bc2'},
 so that it starts on $\tilde{E}^\u$ near the origin at $x_-$
 and arrives on $\tilde{E}^\s$ near the origin at $x_+$,
 where $x_\pm$ are the same as in the above.
So we look for a solution to \eqref{eqn:eig_eq_2'} satisfying
\begin{align}
\tilde{L}^\s
\begin{pmatrix}
\zeta_\R(x_-)\\
		\zeta_\I(x_-)
	\end{pmatrix}
	=0,\quad
	\tilde{L}^\u
	\begin{pmatrix}
		\zeta_\R(x_+)\\
		\zeta_\I(x_+)
	\end{pmatrix}
	=0,
	\label{eqn:bc2}
\end{align}
where
\[
\tilde{L}^\s\defeq
\left(\begin{array}{c|c|c|c}
\begin{array}{cc}
-R_+ & -\Delta_+\\
\Delta_+ & -R_+
\end{array}
& I_4 &
\begin{array}{cc}
\Delta_+ & -R_+\\
R_+ & \Delta_+
\end{array}
& J_4\\\hline
\begin{array}{cc}
-R_- & \Delta_-\\
-\Delta_- & -R_-
\end{array}
& I_4 &
\begin{array}{cc}
\Delta_- & R_-\\
 -R_- & \Delta_-
\end{array}
& -J_4
\end{array}\right),
\]
and
\[
\tilde{L}^\u\defeq
\left(\begin{array}{c|c|c|c}
\begin{array}{cc}
R_+ & \Delta_+\\
-\Delta_+ & R_+
\end{array}
& I_4 &
\begin{array}{cc}
-\Delta_+ & R_+\\
-R_+ & -\Delta_+
\end{array}
& J_4\\\hline
\begin{array}{cc}
R_- & -\Delta_-\\
\Delta_- & R_-
\end{array}
& I_4 &
\begin{array}{cc}
-\Delta_- & -R_-\\
R_- & -\Delta_-
\end{array}
& -J_4
\end{array}\right)
\]
with
\[
R_\pm\defeq
\begin{pmatrix}
\rho_{1\pm} & 0\\[1ex]
0 & \rho_{2\pm}
\end{pmatrix},\quad
\Delta_\pm\defeq
\begin{pmatrix}
\delta_{1\pm} & 0\\[1ex]
0 & \delta_{2\pm}
\end{pmatrix},\quad
J_4\defeq
\begin{pmatrix}
O_2 & I_2\\[1ex]
-I_2 & O_2
\end{pmatrix}.
\]
Here $\tilde{L}^\s$ and $\tilde{L}^\u$ are $8\times16$ real matrices consisting of
 bases in the subspaces spanned by row eigenvectors for the $16\times16$ matrix
\begin{align*}
	\begin{pmatrix}
		A_{\R\infty}(\lambda_\R,\lambda_\I) & -A_{\I\infty}(\lambda_\R,\lambda_\I)\\
		A_{\I\infty}(\lambda_\R,\lambda_\I) & A_{\R\infty}(\lambda_\R,\lambda_\I)
	\end{pmatrix}
\end{align*}
such that the associated eigenvalues have negative and positive real parts, respectively.
We have also denoted
\begin{align*}
	\sqrt{\omega\pm\i\lambda}
	=\rho_{1\pm}+\i\delta_{1\pm},\quad
	\sqrt{s\pm\i\lambda}
	=\rho_{2\pm}+\i\delta_{2\pm}
\end{align*}
with
\begin{align*}
	&\rho_{1\pm}=\sqrt{\frac{\sqrt{(\omega\mp\lambda_\I)^2+\lambda_\R^2}+(\omega\mp\lambda_\I)}{2}},\\
	&\delta_{1\pm}=\pm\sgn(\lambda_\R)\sqrt{\frac{\sqrt{(\omega\mp\lambda_\I)^2+\lambda_\R^2}-(\omega\mp\lambda_\I)}{2}},\\
	&\rho_{2\pm}=\sqrt{\frac{\sqrt{(s\mp\lambda_\I)^2+\lambda_\R^2}+(s\mp\lambda_\I)}{2}},\\
	&\delta_{2\pm}=\pm\sgn(\lambda_\R)\sqrt{\frac{\sqrt{(s\mp\lambda_\I)^2+\lambda_\R^2}-(s\mp\lambda_\I)}{2}}.
\end{align*}
Unlike $|z(x_\pm)|$,
 the distances $|\zeta_{\R}(x_\pm)|,|\zeta_{\I}(x_\pm)|$ do not have to be kept small
 necessarily in the computation since if $(\zeta_{\R}(x_+),\zeta_{\I}(x_+))\in\tilde{E}^\s$
 (resp.\ $(\zeta_{\R}(x_-),\zeta_{\I}(x_-))\in\tilde{E}^\u$) for $|x_\pm|\gg 1$,
 then $(\zeta_{\R}(x),\zeta_{\I}(x))$ tends to the origin
 as $x\to+\infty$ (resp.\ $x\to-\infty$).
To eliminate the multiplicity of solutions due to the linearity of \eqref{eqn:eig_eq},
 we also add the integral conditions
\begin{align}
	\begin{aligned}
		&\sum_{j=1}^4\int_{x_-}^{x_+}\bigl((\zeta_{\R j}(x)-\zeta^\ast_{\R j}(x))\zeta^\ast_{\R j}(x)
			+(\zeta_{\I j}(x)-\zeta^\ast_{\I j}(x))\zeta^\ast_{\I j}(x)\bigr)\d x=0,\\
		&\sum_{j=1}^4\int_{x_-}^{x_+}\bigl((\zeta_{\I j}(x)-\zeta^\ast_{\I j}(x))\zeta^\ast_{\R j}(x)
			-(\zeta_{\R j}(x)-\zeta^\ast_{\R j}(x))\zeta^\ast_{\I j}(x)\bigr)\d x=0,
	\end{aligned}
	\label{eqn:ic2}
\end{align}
which are equivalent to
\begin{align*}
	\sum_{j=1}^4\int_{x_-}^{x_+}\bigl(\zeta_j(x)-\zeta^*_j(x)\bigr)\overline{\zeta^*_j(x)}\,\d x=0,
\end{align*}
where $\zeta^*=(\zeta^*_1,\ldots,\zeta^*_8)^\T$,
 $\zeta_\R^*=(\zeta^*_{\R1},\ldots,\zeta^*_{\R 8})^\T$
 and $\zeta_\I^*=(\zeta^*_{\I1},\ldots,\zeta^*_{\I 8})^\T$ represent
 previously computed solutions along continuation branches.


\subsection{Numerical Results}

We used the computer continuation tool \texttt{AUTO} \cite{D12}
 to obtain numerical solutions to \eqref{eqn:dyn_sys} and \eqref{eqn:eig_eq_2'}
 satisfying the boundary conditions \eqref{eqn:bc1} and \eqref{eqn:bc2}, respectively,
 under the integral conditions \eqref{eqn:ic1} and \eqref{eqn:ic2}, as in \cite{YY20a}.
In the numerical continuations, $\beta_1$ was varied
 along with $d_1$, $\lambda_\R$ and $\lambda_\I$ taken as free parameters.
Moreover, the homoclinic solution \eqref{eqn:homo}
 and the eigenfunctions \eqref{eqn:ef1} or \eqref{eqn:ef2}
 with the eigenvalues \eqref{eqn:ev} were taken as a starting solution.
The distances $|z(x_\pm)|$ were monitored 
 and kept small ($\approx10^{-3}$ typically).

\begin{table}
\caption{Values of $\beta_1^{(\ell)}$ and $\bar{b}_2$
 for $(\omega,s,\beta_2)=(1,4,2)$.
The values of $\bar{b}_2$ are rounded off to the fourth place.}
\label{tab:ab}
\renewcommand{\arraystretch}{1.4}
\begin{tabular}{|c|c|c|c|c|c|} \hline
$\ell$ & \makebox[4em]{0} & \makebox[4em]{1} & \makebox[4em]{2}
 & \makebox[4em]{3} & \makebox[4em]{4} \\\hline
$\beta_1^{(\ell)}$ & 3 & 6 & 10 & 15 & 21\\\hline
\raisebox{-1pt}{$\bar{b}_2$} & $5.486$ & $0.3879$ & $0.03650$ & $0.001333$ & $-0.002094$ \\\hline
\end{tabular}
\end{table}

We set $\omega=1$, $s=4$ and $\beta_2=2$.
The constant $\bar{b}_2$ appearing in Theorem~\ref{thm:bif} and $\beta_1^{(\ell)}$
 were calculated according to the formulas given in Appendix~B of \cite{YY20b} and \eqref{eqn:beta1_ell}
 as in Table~\ref{tab:ab}.
From Theorem~\ref{thm:bif} and Table~\ref{tab:ab}
 we see that the first four pitchfork bifurcations are supercritical
 but the fifth one is subcritical.

\begin{figure}[t]
	\begin{minipage}{0.495\textwidth}
		\includegraphics[scale=0.55]{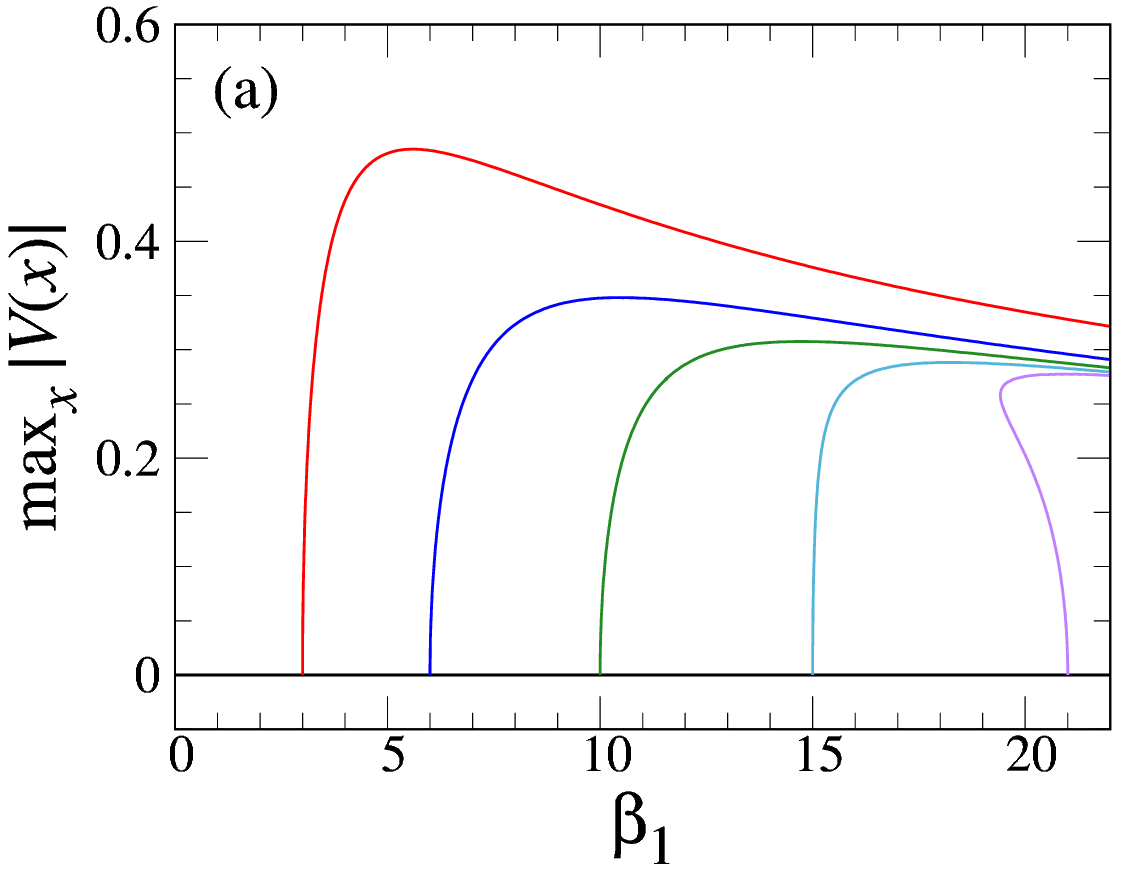}
	\end{minipage}
	\begin{minipage}{0.495\textwidth}
	\begin{center}
		\includegraphics[scale=0.45]{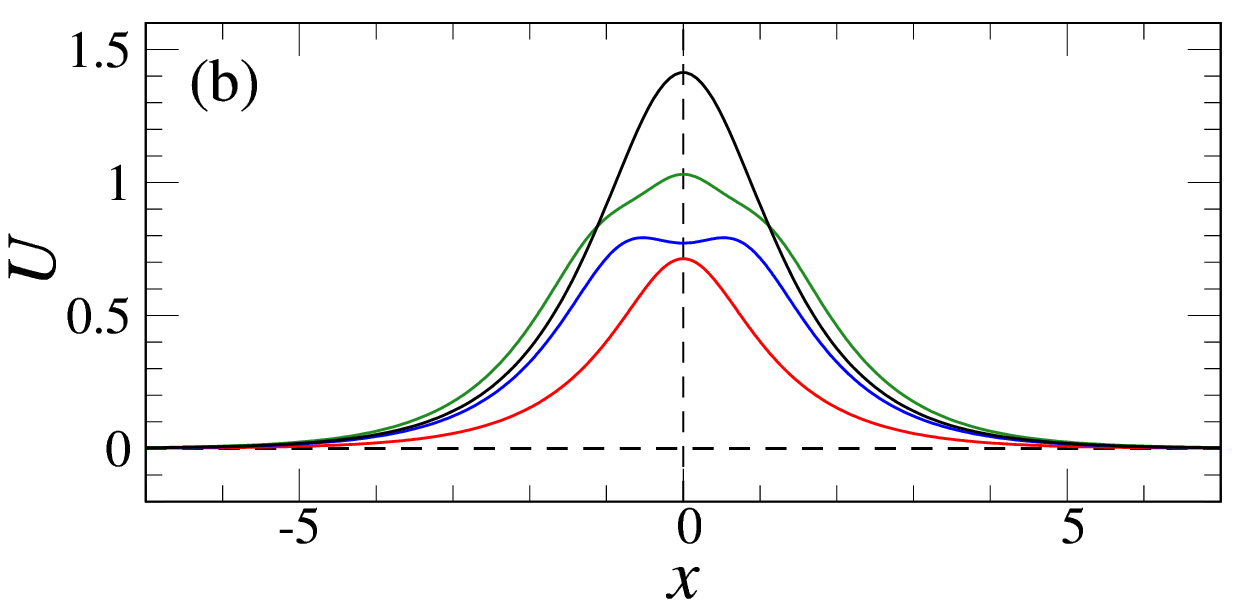}\\[2ex]
		\includegraphics[scale=0.45]{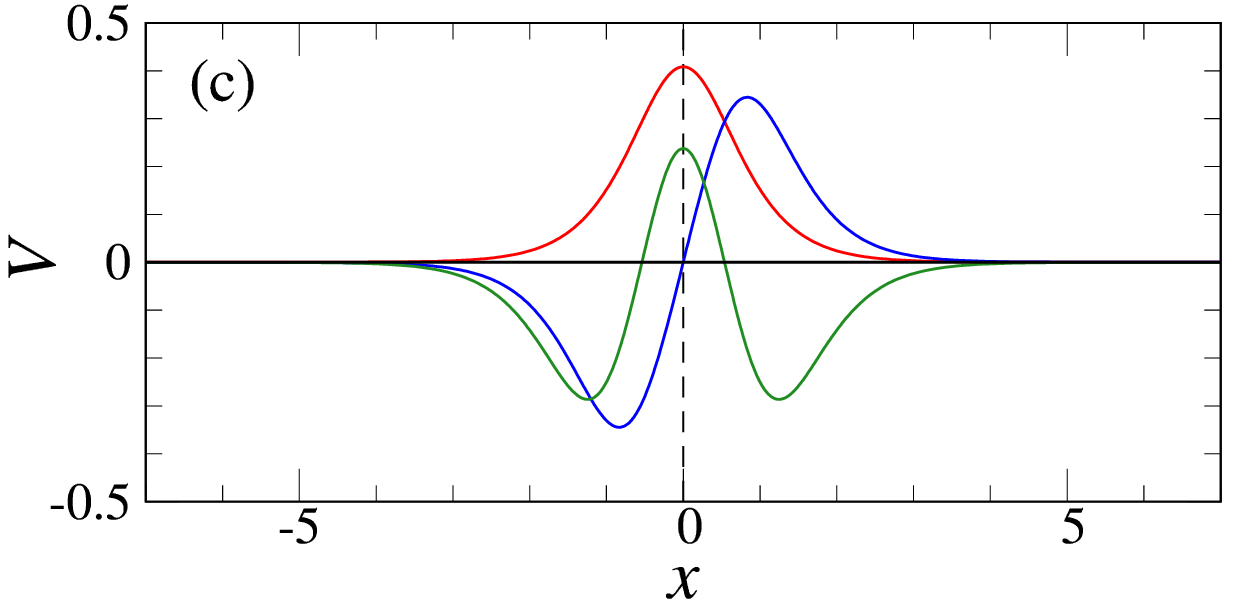}
	\end{center}
	\end{minipage}
	\begin{center}
		\includegraphics[scale=0.45]{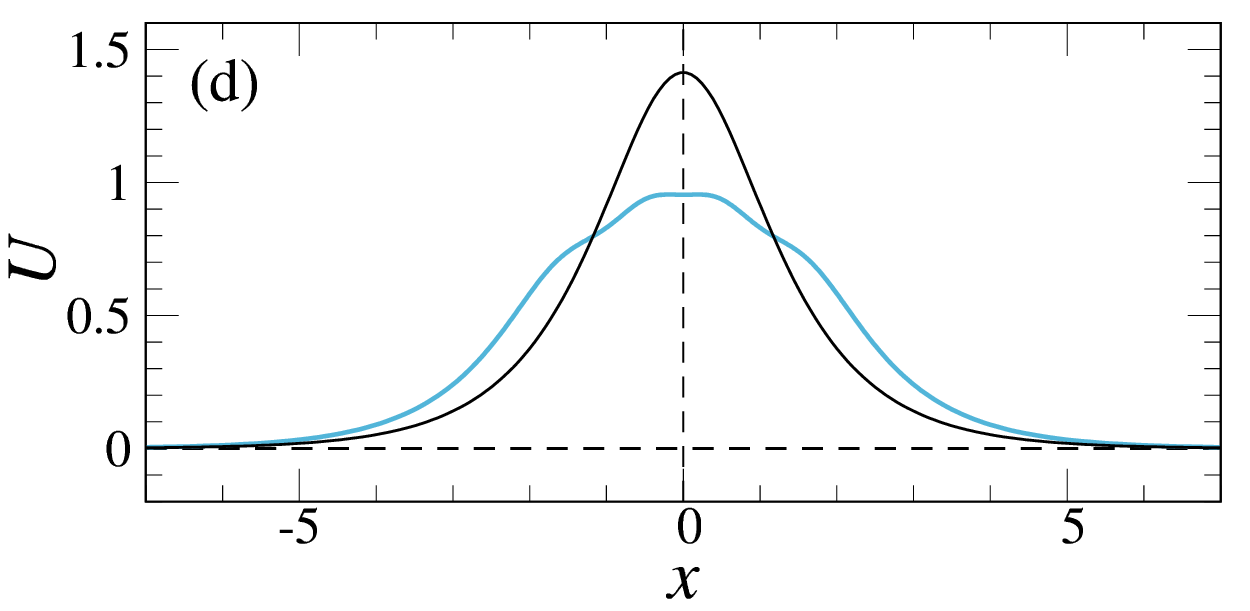}\quad
		\includegraphics[scale=0.45]{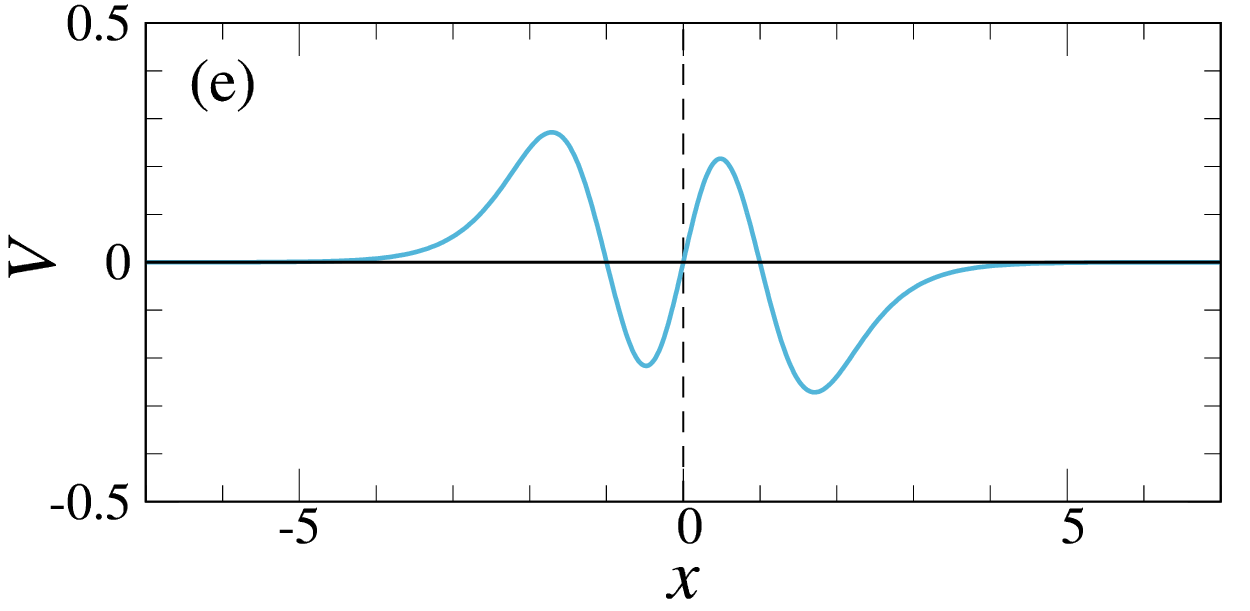}
	\end{center}
	\caption{
		Bifurcations of solitary waves in \eqref{eqn:CNLS}
		for $(\omega,s,\beta_2)=(1,4,2)$:
		(a) Bifurcation diagram;
		(b) and (c) (resp.\ (d) and (e)): profiles of the corresponding homoclinic solutions
		(resp.\ solution) to \eqref{eqn:SE} on the first three branches at $\beta_1=12$
		(resp.\ on the fourth branch at $\beta_1=16$).
		In plate~(a), the red, blue, green, light blue, and purple lines
		represent the branches born at the first, second, third, fourth and fifth
		bifurcations (at $\beta_1=3,6,10,15$, and $21$), respectively,
		while the black line represents the branch of the fundamental solitary wave
		\eqref{eqn:fun_sw}.
		In plates~(b)-(e), the homoclinic solutions along with \eqref{eqn:homo}
		are plotted as the same color lines as the corresponding branches in plate~(a).
		See Fig.~\ref{fig:gker_UV} for profiles of the corresponding
		homoclinic solution to \eqref{eqn:SE} on the fifth branch.}
	\label{fig:bif_diag}
\end{figure}

\begin{figure}[thb]
	\begin{center}
		\includegraphics[scale=0.49]{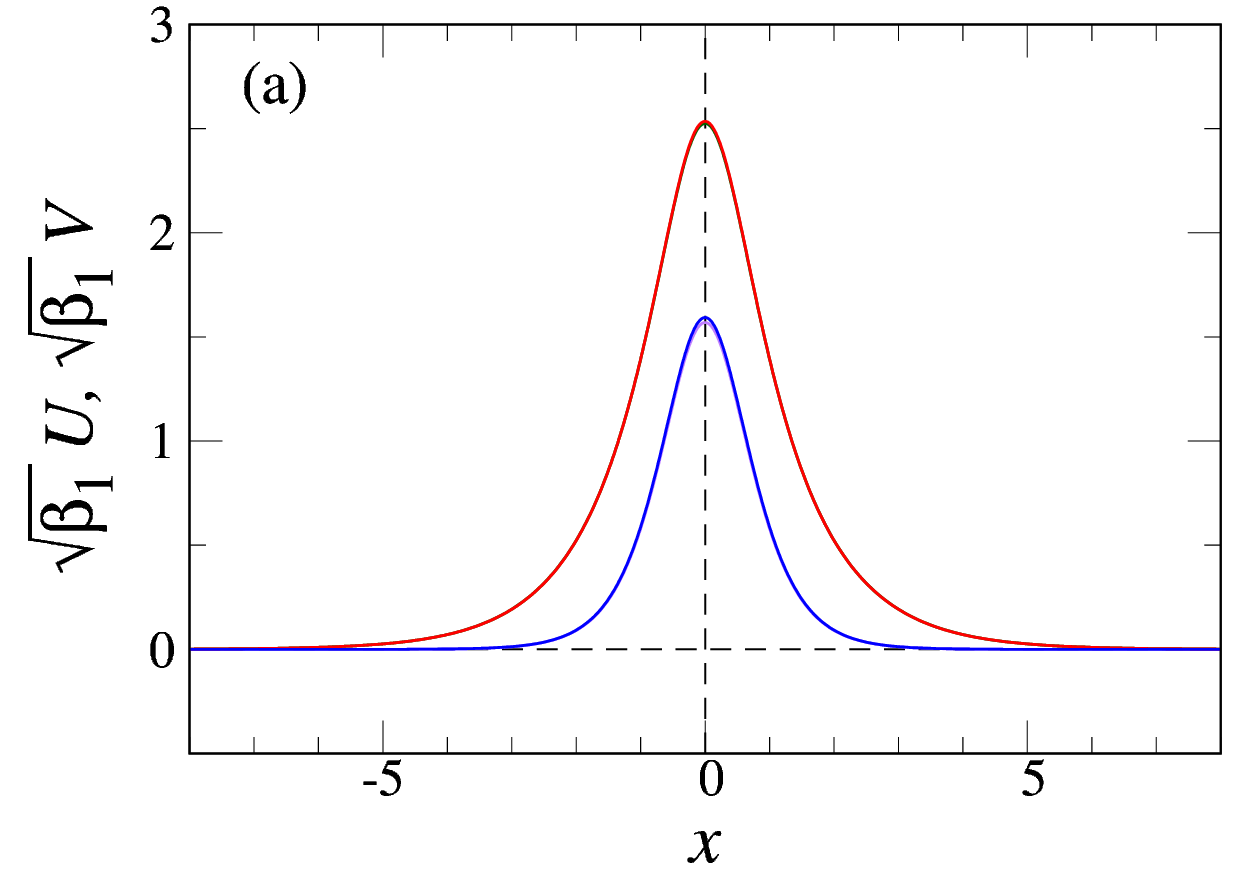}
		\includegraphics[scale=0.49]{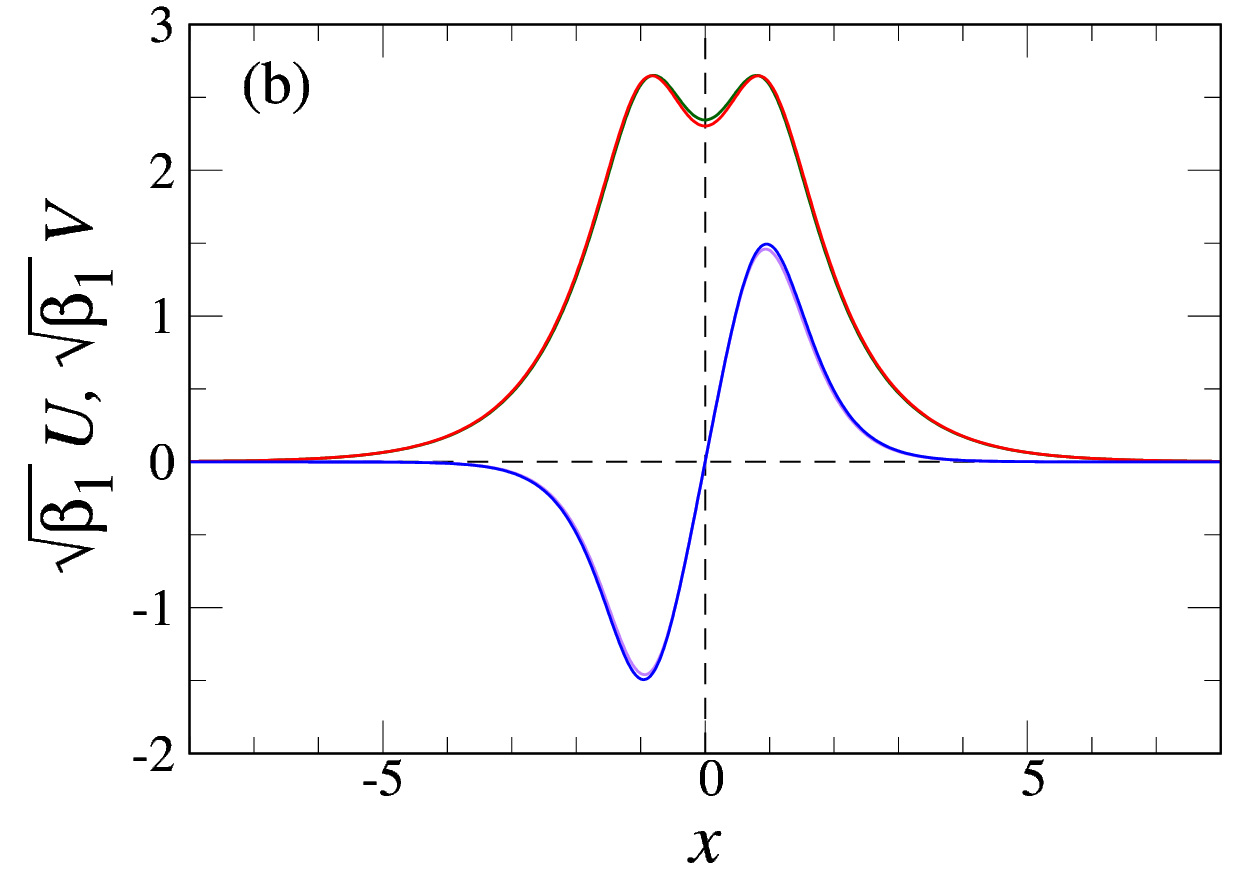}\\[1ex]
		\includegraphics[scale=0.49]{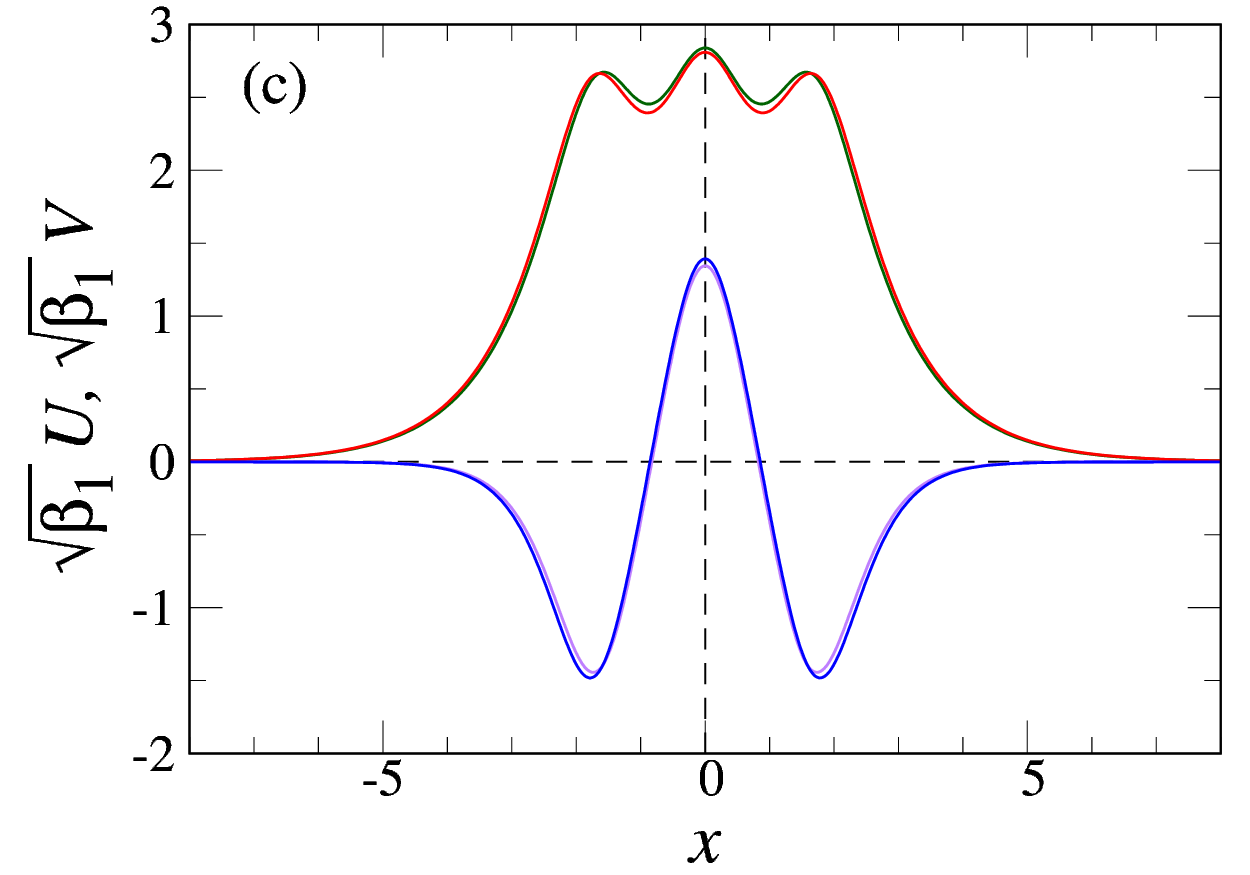}
		\includegraphics[scale=0.49]{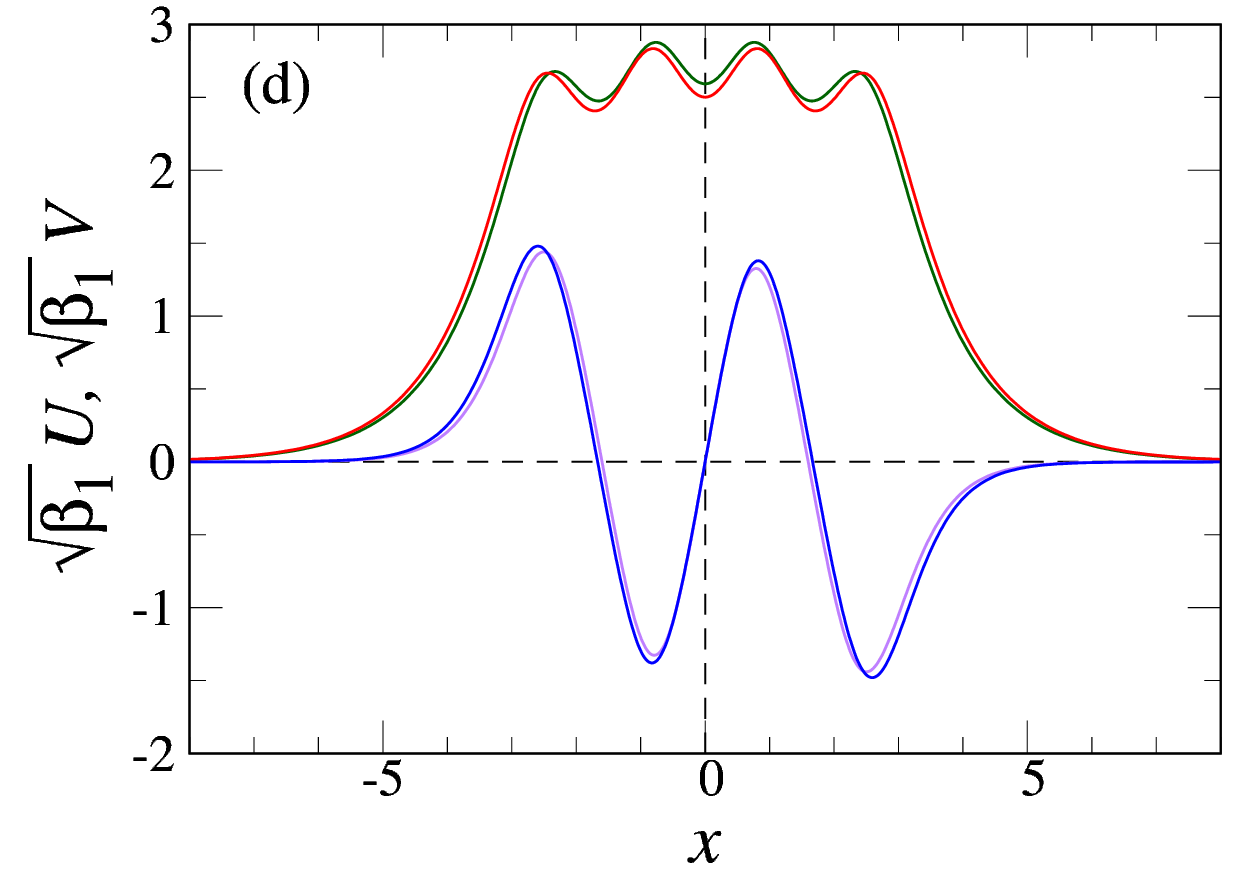}\\[1ex]
		\includegraphics[scale=0.49]{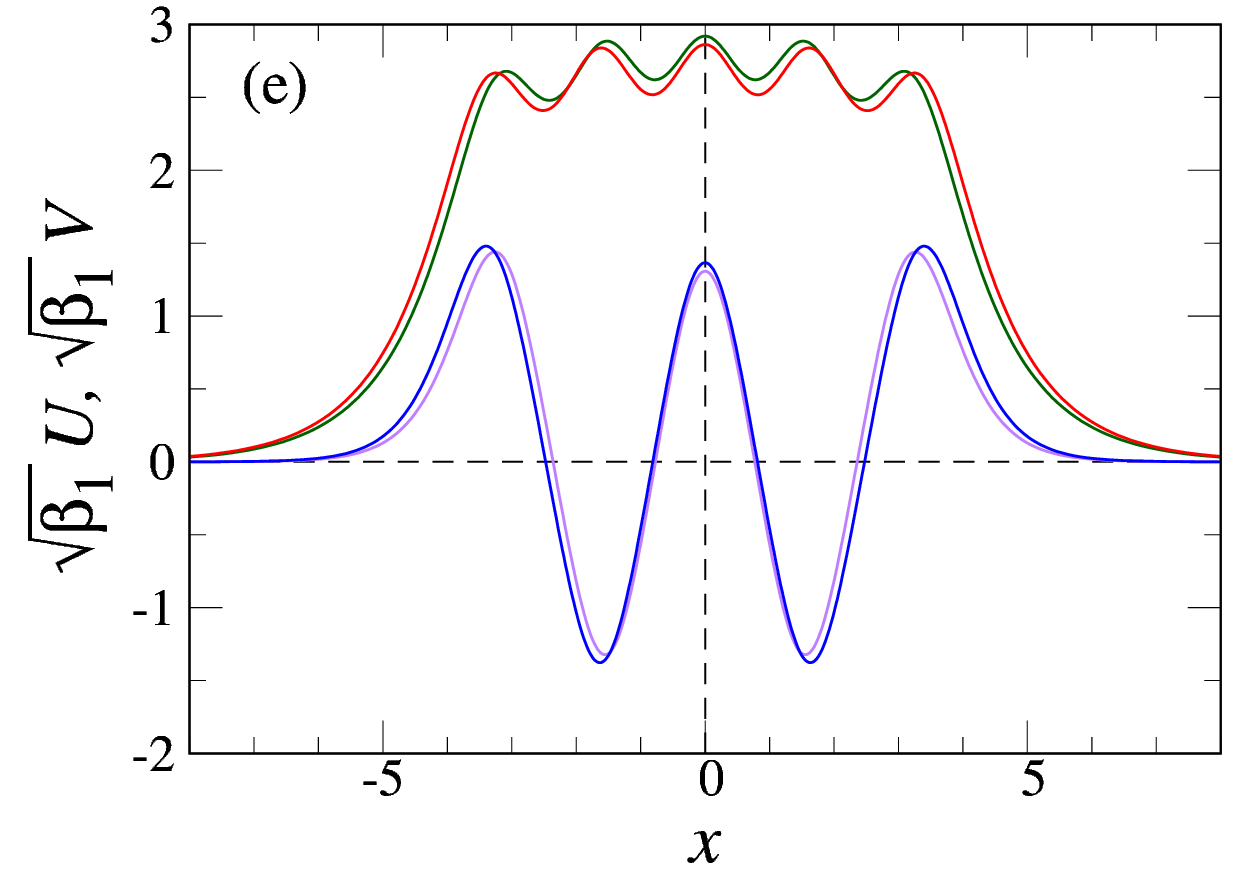}
	\end{center}
	\caption{
		Profiles of the homoclinic solutions to \eqref{eqn:SE} on each branch
		at $\beta_1=50$ and $100$:
		(a) $\ell=0$; (b) $\ell=1$; (c) $\ell=2$; (d) $\ell=3$; (e) $\ell=4$.
		The red and blue (resp.\ green and purple) lines, respectively, represents
		the $U$- and $V$-components at $\beta_1=100$ (resp.\ at $\beta_1=50$).
	}
	\label{fig:asym}
\end{figure}

Figure \ref{fig:bif_diag}(a) shows a numerically computed bifurcation diagram
 of homoclinic solutions to \eqref{eqn:SE},
 which correspond to  solitary waves in \eqref{eqn:CNLS},
 where $x_\pm=\pm7$ were taken
 except that $x_\pm=\pm5$ for the first branch ($\ell=0$)
 because of rapid dcaying of the $V$-component as $x\to\pm\infty$ .
We observe that pitchfork bifurcations of solitary waves occur
 at $\beta_1=3,6,10,15,21$, as predicted in Theorem \ref{thm:bif}
 (see also Table~\ref{tab:ab}).
Note that a pair of symmetric branches about $V=0$ are born at each bifurcation point.
Moreover, we see that a saddle-node bifurcation of the solitary waves
 on the fifth branches occurs at $\beta_1\approx 19.4$
 (more precisely, $19.41626\ldots$).
The homoclinic solutions to \eqref{eqn:SE} on the branches
 born at the first four bifurcation points 
 are displayed for $\beta_1=12$ and $\beta_1=16$ 
 along with the homoclinic solution \eqref{eqn:homo}
 in Figs.~\ref{fig:bif_diag}(b)-(e).
The $V$-component of the homoclinic orbit on the $(\ell+1)$th branch
 have exactly $\ell$ zeros for $\ell=0$-$4$.
 
The profiles of the bifurcated homoclinic solutions
 on each branch at $\beta_1=50$ and $100$ are also plotted
 with a scaling of $1/\sqrt{\beta_1}$ in Figure~\ref{fig:asym}.
Here $x_\pm=\pm8$ were used
 since some homoclinic solutions do not decay in a long interval
 (see Figs.~\ref{fig:asym}(d) and (e)).
Thus, they converge to certain shapes
 with a scaling of $1/\sqrt{\beta_1}$ as $\beta_1\to\infty$. 
 
\begin{figure}[thbp]
	\begin{minipage}{0.495\textwidth}
		\includegraphics[scale=0.55]{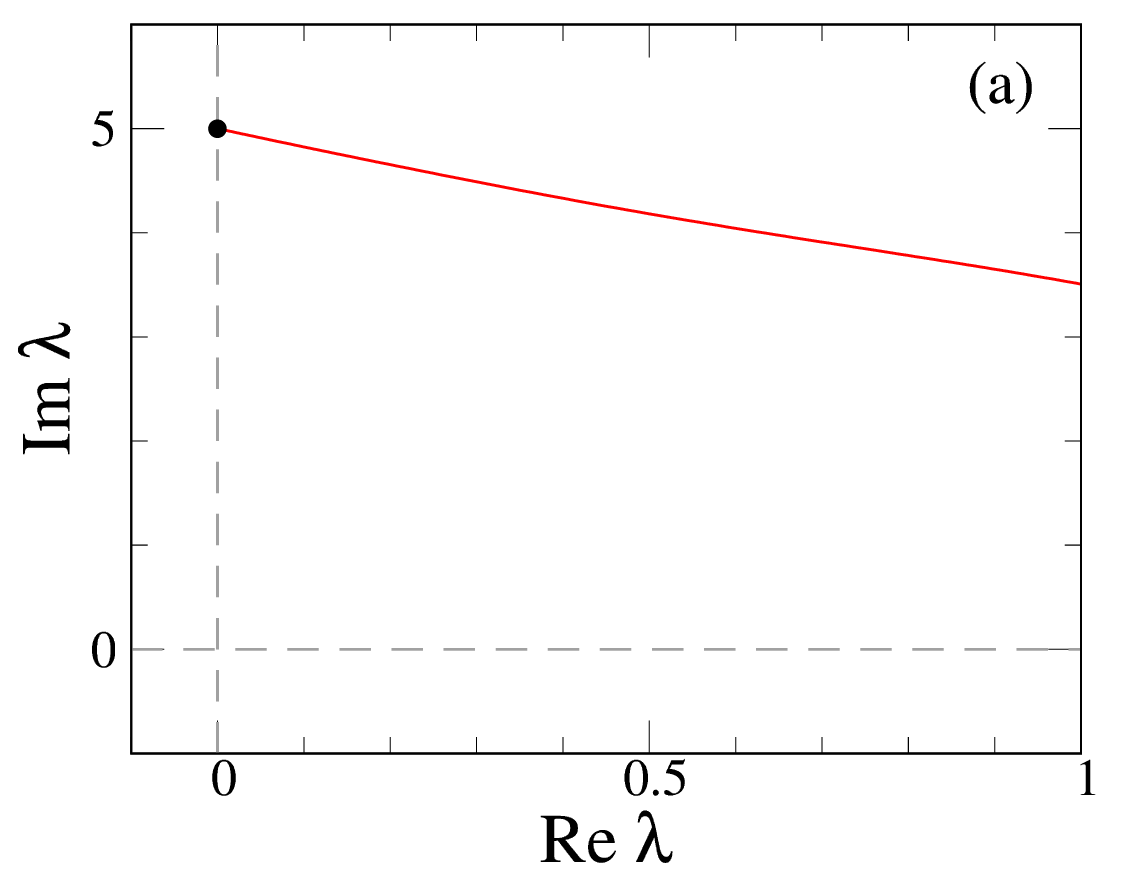}\\[1ex]
		\includegraphics[scale=0.55]{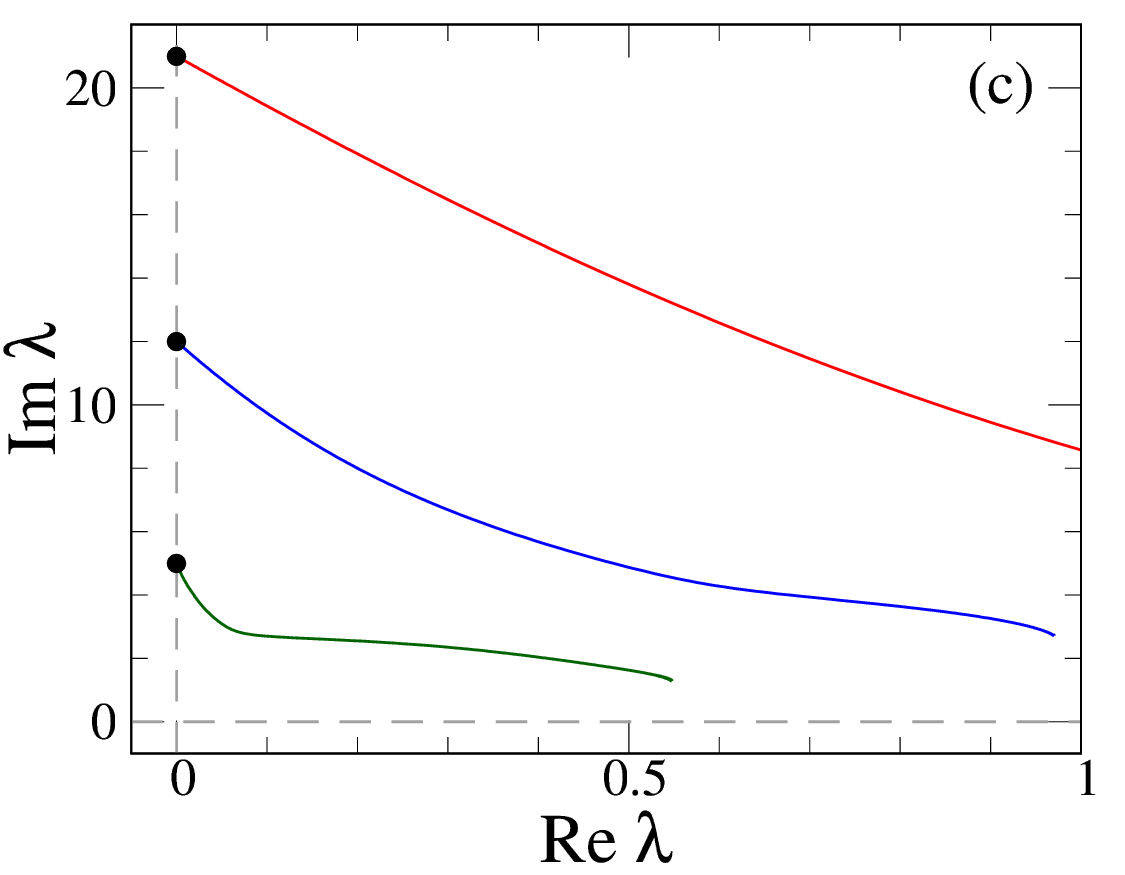}
	\end{minipage}
	\begin{minipage}{0.495\textwidth}
	\begin{center}
		\includegraphics[scale=0.55]{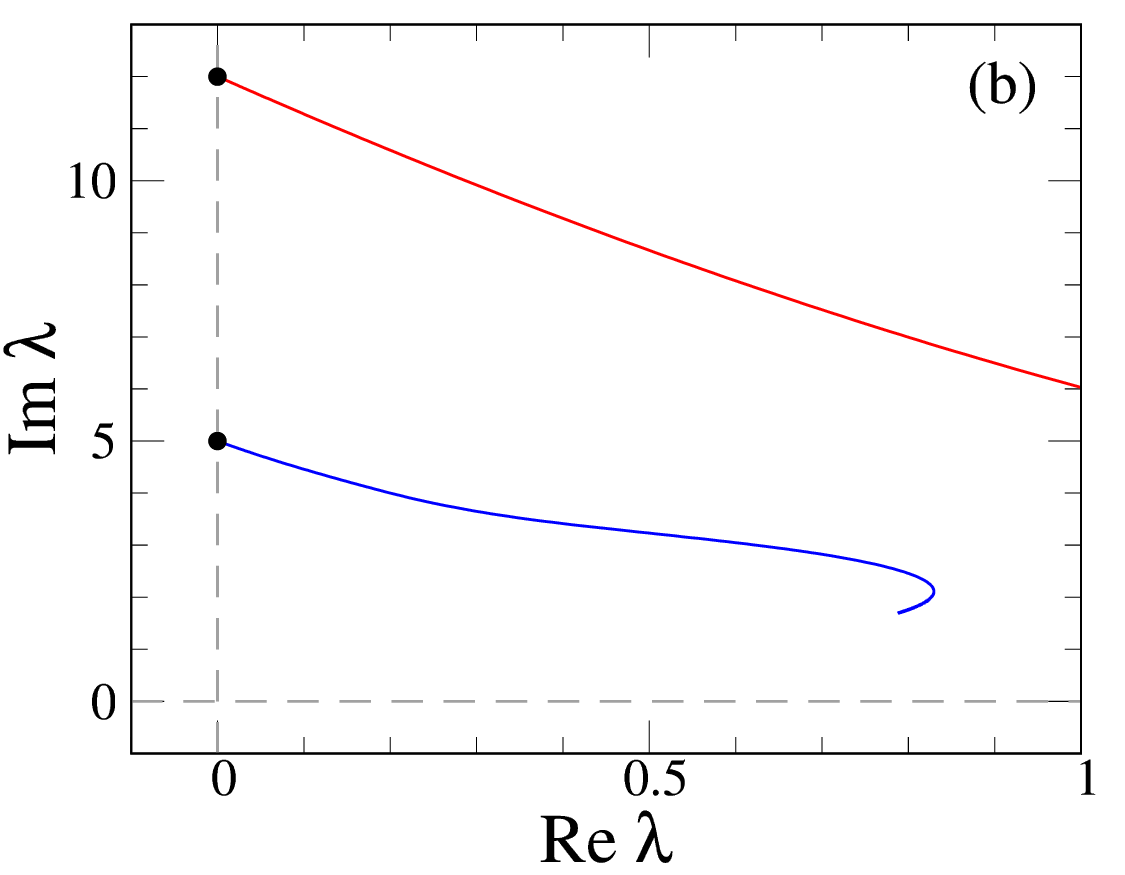}\\[1ex]
		\includegraphics[scale=0.55]{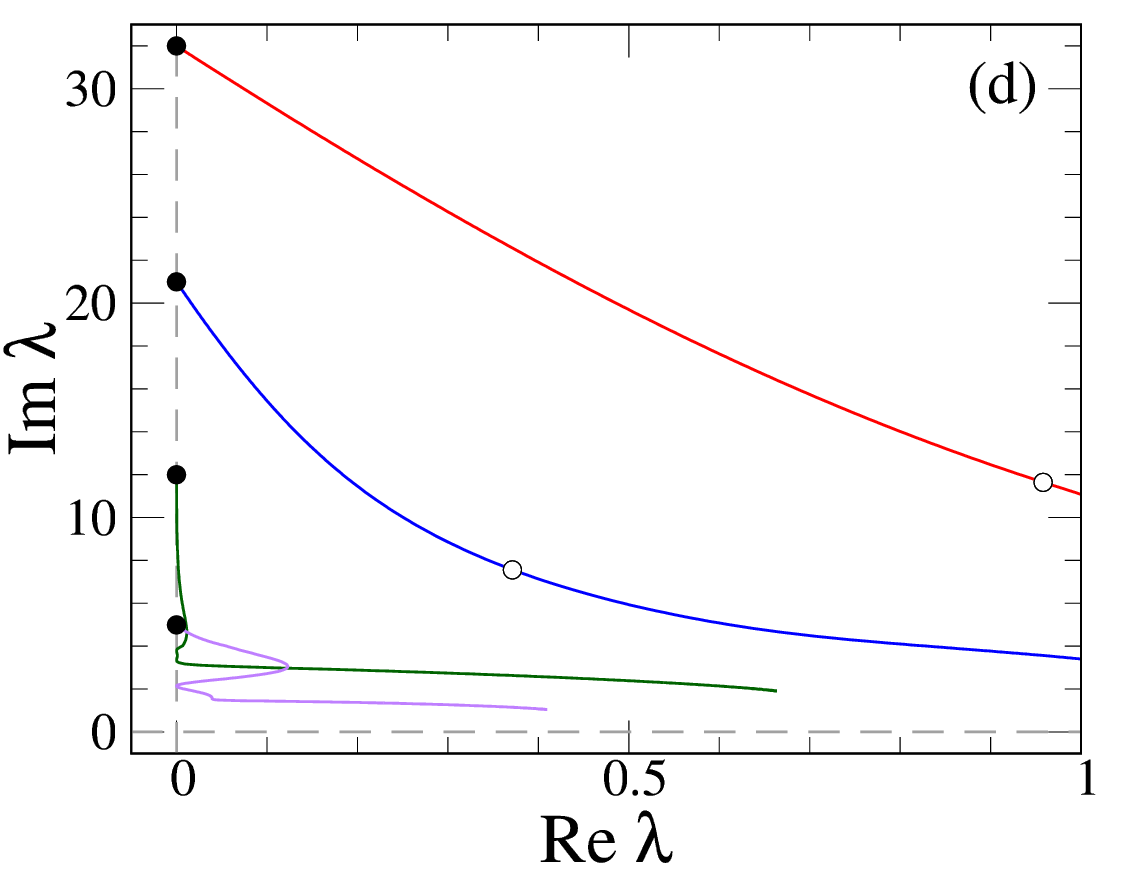}
	\end{center}
	\end{minipage}
	\begin{center}
		\includegraphics[scale=0.55]{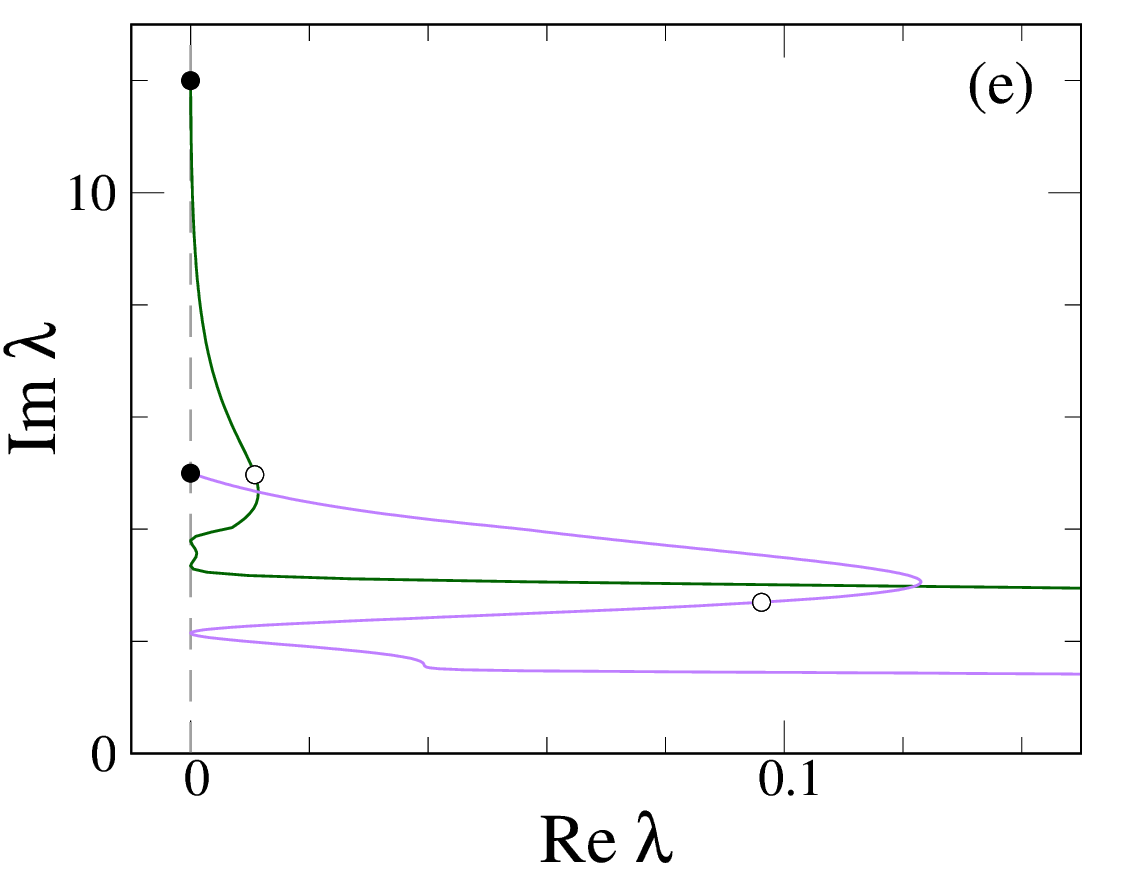}
	\end{center}
	\caption{
		Eigenvalues of the linearized operator $J\mathcal L$ around the bifurcated solitary wave born at $\beta_1=\beta_1^{(\ell)}$ for $(\omega,s,\beta_2)=(1,4,2)$:
		(a) $\ell=1$; (b) $\ell=2$; (c) $\ell=3$; (d) and (e) $\ell=4$.
		Plate~(e) is an enlargement of plate~(d).
		The red, blue, green, and purple lines represent the eigenvalues of $J\mathcal L$
		 with $k=0,1,2$, and $3$, respectively.
		The bullet `$\bullet$' represents the loci of the eigenvalues
		 at $\beta_1=\beta_1^{(\ell)}$.
		In plates (d) and (e), the circle `$\circ$' represents the loci of the eigenvalues
		 at the saddle-node bifurcation point $\beta_1\approx19.41626$.
		Each curve was computed from $\beta_1=\beta_1^{(\ell)}$ to $100$.
	}
	\label{fig:ev}
\end{figure}

\begin{figure}[thb]
	\begin{center}
		\includegraphics[scale=0.49]{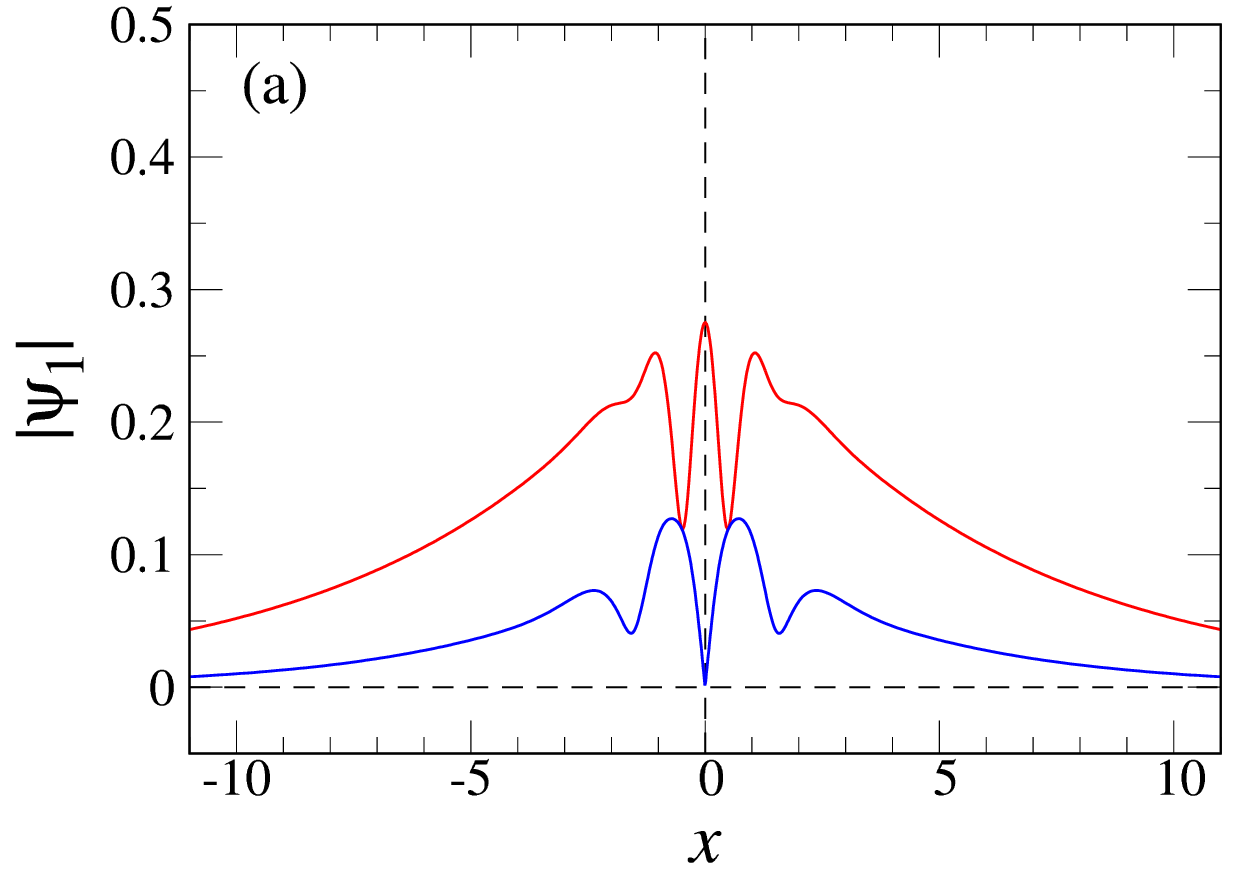}
		\includegraphics[scale=0.49]{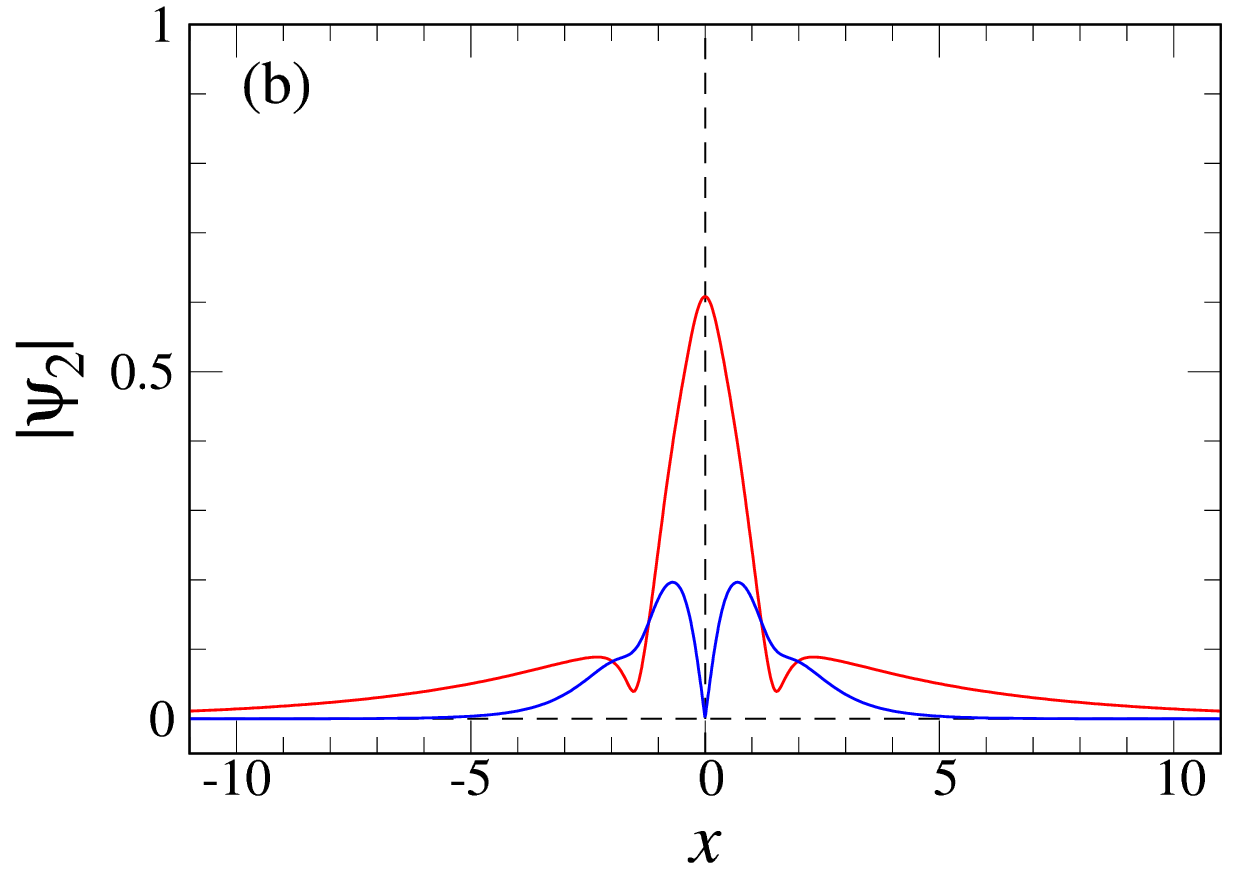}\\[1ex]
		\includegraphics[scale=0.49]{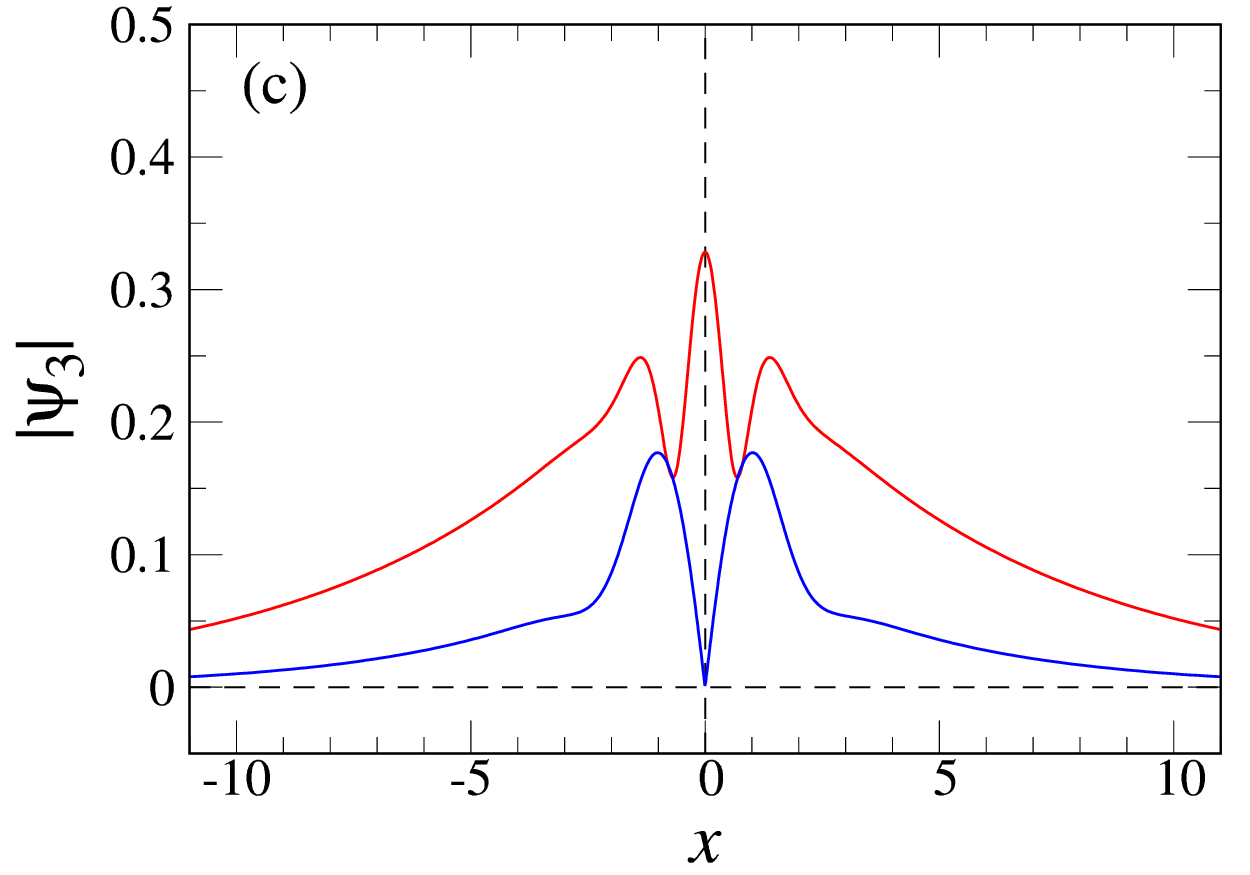}
		\includegraphics[scale=0.49]{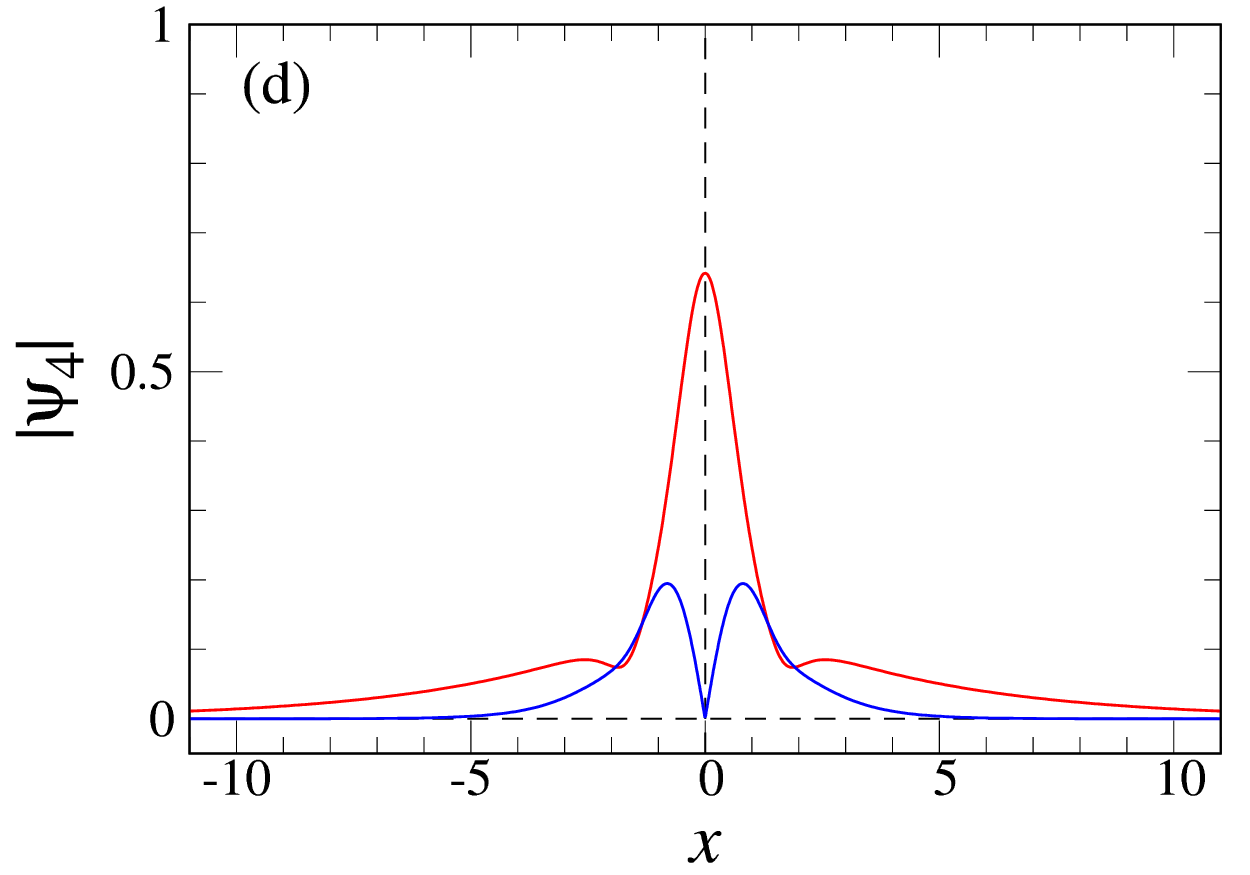}\\[1ex]
		\includegraphics[scale=0.49]{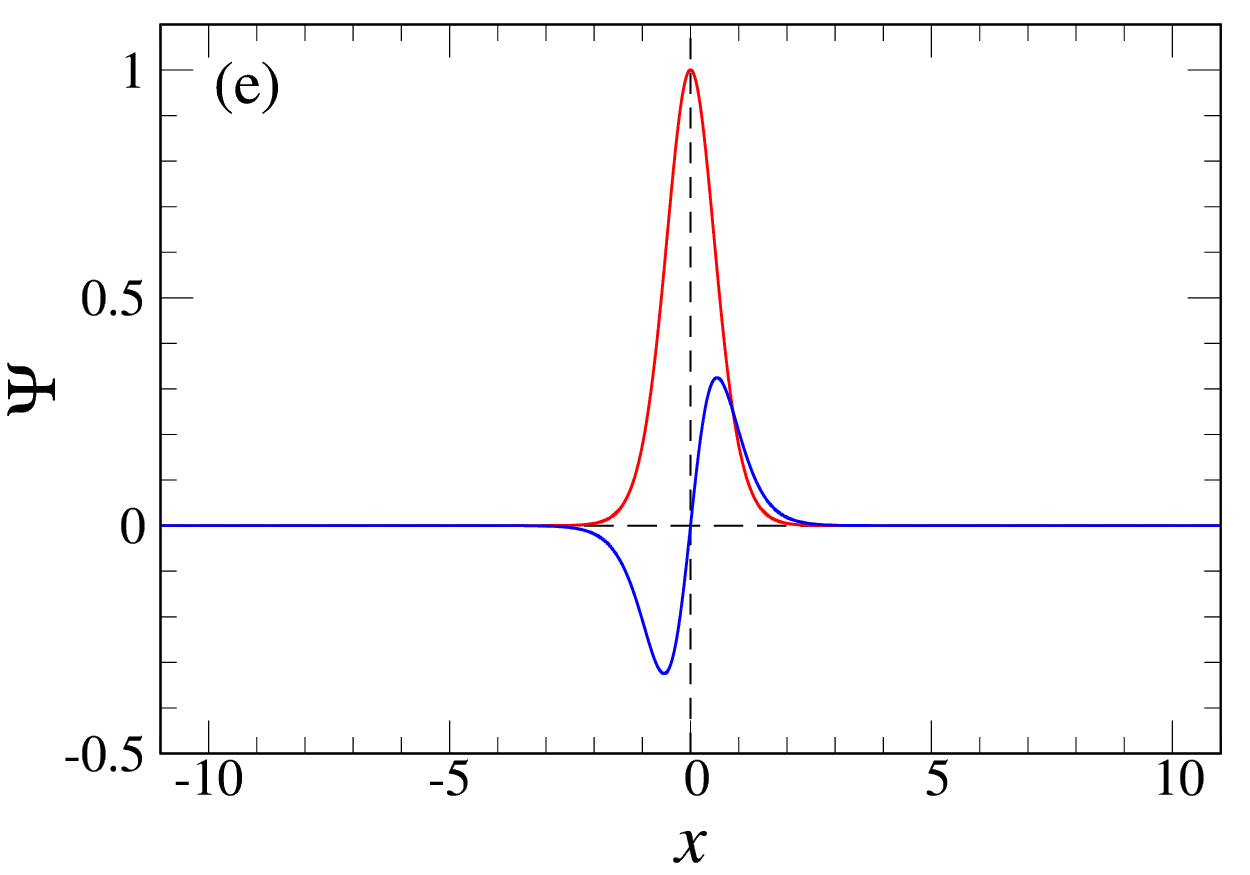}
	\end{center}
	\caption{
		Profiles of the eigenfunctions of the eigenvalue problem \eqref{eqn:eig_eq} for $\ell=2$ at $\beta_1=12$ with $(\omega,s,\beta_2)=(1,4,2)$:
		(a) $|\psi_1|$;  (b) $|\psi_2|$; (c) $|\psi_3|$; (d) $|\psi_4|$; (e) $\Psi$.
		The red and blue lines, respectively, represent the eigenfunctions for the curves emerging from $12\i$ and $5\i$ in Fig.~\ref{fig:ev}(b).
		\label{fig:eig_fun}
	}
\end{figure}

Figure \ref{fig:ev} shows how the eigenvalues of $J\mathcal L$
 for the bifurcated solitary wave on each branch change,
 where $x_\pm=\pm11$ was taken.
We only display the eigenvalues with $\Re\lambda,\Im\lambda\ge0$
 since the spectra of $J\mathcal L$ are symmetric, as stated in Section~2.
For the bifurcated solitary wave born at $\beta_1=\beta_1^{(\ell)}$,
 all eigenvalues of $J\mathcal L$ at $\beta_1=\beta_1^{(\ell)}$ are given by
\begin{align}
\lambda=&\pm\i(s-\omega(\sqrt{s/\omega}+\ell-k)^2),\notag\\
& k\in\{0,1,\ldots,\lfloor\sqrt{s/\omega}\rfloor+\ell\}\setminus\{\sqrt{s/\omega}+\ell\},
\label{eqn:lambda}
\end{align}
(see \eqref{eqn:ev}) and they immediately disappear for $k\ge\ell$
 when $\beta_1$ changes from $\beta_1^{(\ell)}$ (see Remark~7.9 of \cite{YY20b}).
The loci of the eigenvalues of $J\mathcal L$ 
 leaving the imaginary axis from \eqref{eqn:lambda}
 when $\beta_1$ changes from $\beta_1^{(\ell)}$
 are plotted for $k=0$ and $\ell=1$ in Fig.~\ref{fig:ev}(a);
 for $k=0,1$ and $\ell=2$ in Fig.~\ref{fig:ev}(b);
 for $k=0,1,2$ and $\ell=3$ in Fig.~\ref{fig:ev}(c);
 and $k=0,1,2,3$ and $\ell=4$ in Figs.~\ref{fig:ev}(d) and (e).
Each curve was computed from $\beta_1=\beta_1^{(\ell)}$ to $100$.
Although \eqref{eqn:dyn_sys} and \eqref{eqn:eig_eq_2} are highly degenerate
 at the bifurcation point $\beta_1=\beta_1^{(\ell)}$
 since two branches of eigenfunctions are also created there,
 continuation of their solutions by \texttt{AUTO} succeeded from there.
These results indicate that
 the real parts of the eigenvalues become positive
 and the bifurcated solitary waves are unstable
 when $\beta_1\ne\beta_1^{(\ell)}$, as stated in Theorem~\ref{thm:ls}.

On the other hand, the computed eigenvalues at $\beta_1=50$ were almost the same as at $\beta_1=100$ for all computed branches.
So the eigenvalues are thought to converge to certain values as $\beta_1\to\infty$.
See, e.g., the blue line ($k=1$) in Fig.~\ref{fig:ev}(b),
 the blue and green lines ($k=1$ and $2$) in Fig.~\ref{fig:ev}(c),
 and the green and purple lines ($k=2$ and $3$)  in Fig.~\ref{fig:ev}(d).
The reason is that as stated above, as $\beta_1\to\infty$,
 the bifurcated homoclinic solutions converge
 to a certain profile, say $(U_\infty(x),V_\infty(x))$, with a scaling of $1/\sqrt{\beta_1}$
 and the operators $\mathcal L_\pm$ in $J\mathcal L$ with \eqref{eqn:JL} converge to
\begin{align*}
	&\mathcal L_+=\begin{pmatrix}
		-\partial_x^2+\omega-V_\infty^2&-2U_\infty V_\infty\\
		-2U_\infty V_\infty&-\partial_x^2+s-U_\infty^2
	\end{pmatrix},\\
	&\mathcal L_-=\begin{pmatrix}
		-\partial_x^2+\omega-V_\infty^2 &0\\
		0&-\partial_x^2+s-U_\infty^2
	\end{pmatrix}.
\end{align*}
Note that $(U,V)\to 0$ while $\sqrt{\beta_1}(U,V)\to(U_\infty,V_\infty)$.

In Figs.~\ref{fig:ev}(d) and (e)
 four eigenvalues for the bifurcated solitary wave on the fifth branch ($\ell=4$)
 are displayed and their values at $\beta_1=\beta_1^\mathrm{SN}\approx 19.41626$,
 at which a saddle-node bifurcation occurs (see Fig.~\ref{fig:bif_diag}(a)),
 are plotted as a circle `$\circ$'.
In particular, the solitary wave seems not to change its stability type
 at the saddle-node bifurcation point
 since all the eigenvalues are far from the imaginary axis.
On the other hand, according to Theorem~2.4 of \cite{BY12},
 the VE \eqref{eqn:VE} around the corresponding homoclinic solution
 has two linearly independent solutions there since no bifurcation occurs if it does not,
 so that the geometrical multiplicity of the zero eigenvalue of $J\mathcal L$
 increases by one.
So we suspect that for the zero eigenvalue
 a generalized eigenfuntion turns to an eigenfunction there,
 as suggested from Yang's result \cite{Ya12b}
 for general single NLS equations with external potentials.
This suspicion will be numerically proven true in the next section.

Figure~\ref{fig:eig_fun} displays the absolute value of each component
 of the corresponding eigenfunction $\psi=(\psi_1,\ldots,\psi_4)^\T$
 for $\ell=2$ at $\beta_1=14$.
In Fig.~\ref{fig:eig_fun}(e), the profiles of $\Psi(x)$ for $\ell=2$,
 which are given by \eqref{eqn:ef1} and \eqref{eqn:ef2} with $k=0$ and $1$,
 and represent the eigenfunctions associated with the eigenvalues
 $\lambda=12\i$ and $5\i$ there
 as $\psi=(0,\Psi,0,\i\Psi)^\T$, are plotted as the red and blue lines, respectively.
We see that the eigenfunctions considerably changes
 from those at the bifurcation point $\beta_1=\beta_1^{(2)}$.


\section{Computations of generalized eigenfunctions for the zero eigenvalues}

In this section, for the eigenvalue problem \eqref{eqn:eig_eq},
 we give some numerical computation results
 for its generalized eigenfunction associated with the zero eigenvalue
 which turns to an eigenfunction at the saddle-node bifurcation point.
These results will demonstrate the correctness of our suspicion stated above
 on the saddle-node bifurcation observed in Fig.~\ref{fig:bif_diag}(a):
 the geometric multiplicity of the zero eigenvalue increases by one
 but the number of eigenvalues counted with their multiplicity
 on the imaginary axis does not change, 
 so that the solitary wave does not change its stability type at the bifurcation point.


\subsection{Numerical Approach}

We first recall that $\varphi_j$ and $\chi_j$, $j=1,2,3$,
 can be the eigenfunctions and generalized eigenfunctions of \eqref{eqn:eig_eq}
 associated with the zero eigenvalue
 which are given by \eqref{eqn:kernel} and \eqref{eqn:def_chi}, respectively, such that $J\mathcal L\chi_j=\varphi_j$.
Moreover, the first and second components of $\chi_1$ are zero
 while the third and fourth components of $\chi_2$ and $\chi_3$ are zero.
Hence, if a generalized eigenfunction $\chi$ for the zero eigenvalue
 becomes an eigenfunction at the saddle-node bifurcation point $\beta_1^\mathrm{SN}$, 
 then $\chi$ can be written as a linear combination of $\chi_2$ and $\chi_3$.
This means that the geometric multiplicity is four at $\beta_1=\beta_1^\mathrm{SN}$ but three at $\beta_1\ne\beta_1^\mathrm{SN}$ while the algebraic multiplicity is six for both cases.
So we want to compute a generalized eigenfunction satisfying
\begin{align}
	J\mathcal L\psi=\epsilon_1((1-\epsilon_2)\varphi_2+\epsilon_2\varphi_3)
	\label{eqn:gker}
\end{align}
where $\epsilon_1,\epsilon_2\in\Rset$ are constants.
If $\epsilon_1=0$, then $\psi$ becomes an eigenfunction.
If $\epsilon_2=0$ and $1$ with $\epsilon_1=1$, then $\psi=\chi_2$ and $\chi_3$, respectively, satisfy \eqref{eqn:gker}.

Combining \eqref{eqn:dyn_sys} and \eqref{eqn:gker}, we write our problem as
\begin{align}
	\begin{pmatrix}
		z'\\
		\eta'
	\end{pmatrix}
	=g(z,\eta;\beta_1,\epsilon_1,\epsilon_2),\quad
	z=(z_1,z_2,z_3,z_4)^\T,\eta=(\eta_1,\eta_2,\eta_3,\eta_4)^\T\in\Rset^4,
	\label{eqn:gker2}
\end{align}
and numerically compute a homoclinic solution to \eqref{eqn:gker2} satisfying
\begin{align}
	\lim_{x\to\pm\infty}z(x)=\lim_{x\to\pm\infty}\eta(x)=0,
	\label{eqn:eta_decay}
\end{align}
where $\eta_j=\psi_j$, $j=1,2$, and
\begin{align*}
&
	g(z,\eta;\beta_1,\epsilon_1,\epsilon_2)\notag\\
&
	\defeq\begin{pmatrix}
		z_3\\
		z_4\\
		\omega z_1-(z_1^2+\beta_1z_2^2)z_1+d_1z_3\\
		sz_2-(\beta_1z_1^2+\beta_2z_2^2)z_2+d_1z_4\\
		\eta_3\\
		\eta_4\\
		\omega\eta_1-(3z_1^2+\beta_1z_2^2)\eta_1-2\beta_1z_1z_2\eta_2
		+\epsilon_1(1-\epsilon_2)z_1+d_2z_3\\
		s\eta_2-2\beta_1z_1z_2\eta_1-(\beta_1z_1^2+3\beta_2z_2^2)\eta_2
		+\epsilon_1\epsilon_2z_2+d_2z_4
	\end{pmatrix}
\end{align*}
with dummy parameters $d_1,d_2$.
Here the first and second components of \eqref{eqn:gker}
 and the third and fourth components of $\psi$ have been eliminated
 since the first and second components of both $\varphi_2$ and $\varphi_3$
 and the third and fourth components of both $\chi_2$ and $\chi_3$ are zero.
The fact that the third and fourth components of $\varphi_2$ are $(U(x),0)$
 and those of $\varphi_3$ are $(0,V(x))$ has also been used (see \eqref{eqn:def_chi}).
We easily see that a homoclinic solution persists in \eqref{eqn:gker2}
 when one of the other parameters changes only if $d_1=d_2=0$,
 as in \eqref{eqn:dyn_sys}.
 
Let $\hat{E}^\s$ and $\hat{E}^\u$ be, respectively,
 the four-dimensional stable and unstable subspaces of the linearized system
 at the origin for \eqref{eqn:gker2},
\begin{align*}
	\begin{pmatrix}
		\delta z'\\
		\delta\eta'
	\end{pmatrix}
	=\D_{z,\eta}g(0,0;\beta_1,\epsilon_1,\epsilon_2)
	\begin{pmatrix}
		\delta z\\
		\delta\eta
	\end{pmatrix}
\end{align*}
with
\begin{align*}
	\D_{z,\eta}g(0,0;\beta_1,\epsilon_1,\epsilon_2)
	=\left(
		\begin{array}{c|c}
		\begin{array}{cccc}
			0&0&1&0\\
			0&0&0&1\\
			\omega&0&d_1&0\\
			0&s&0&d_1
		\end{array} & O_4\\\hline
		\begin{array}{cccc}
			0&0&0&0\\
			0&0&0&0\\
			\epsilon_1(1-\epsilon_2)&0&d_2&0\\
			0&\epsilon_1\epsilon_2&0&d_2
		\end{array} &
		\begin{array}{cccc}
			0&0&1&0\\
			0&0&0&1\\
			\omega&0&0&0\\
			0&s&0&0
		\end{array}
		\end{array}
	\right).
\end{align*}
As in Section~3,
 we approximate the homoclinic solution $(z(x),\eta(x))^\T$
 to \eqref{eqn:gker2} satisfying \eqref{eqn:eta_decay},
 so that it starts on $\hat{E}^\u$ near the origin at $x_-$
 and arrives on $\hat{E}^\s$ near the origin at $x_+$.
So we look for a solution to \eqref{eqn:gker2} satisfying
\begin{align}
	\hat{L}^\s
	\begin{pmatrix}
		z(x_-)\\
		\eta(x_-)
	\end{pmatrix}
	=0,\quad
	\hat{L}^\u
	\begin{pmatrix}
		z(x_+)\\
		\eta(x_+)
	\end{pmatrix}
	=0,
	\label{eqn:bc3}
\end{align}
where
\begin{align*}
	&\hat{L}^\s\defeq
	\left(
		\begin{array}{c|c}
		\begin{array}{cc}
			-d_1/2-\sqrt{\omega+(d_1/2)^2}&0\\
			0&-d_1/2-\sqrt{s+(d_1/2)^2}\\
			\hat{L}^\s_{31}&0\\
			0&\hat{L}^\s_{42}
		\end{array} & \hat{L}_-
		\end{array}
	\right),\\
	&\hat{L}^\u\defeq
	\left(
		\begin{array}{c|c}
		\begin{array}{cc}
			-d_1/2+\sqrt{\omega+(d_1/2)^2}&0\\
			0&-d_1/2+\sqrt{s+(d_1/2)^2}\\
			\hat{L}^\u_{31}&0\\
			0&\hat{L}^\u_{42}
		\end{array} & \hat{L}_+
		\end{array}
	\right)
\end{align*}
with
\begin{equation}
\begin{split}
	&\hat{L}^\s_{31}
	=\frac{-\epsilon_1(1-\epsilon_2)-d_2\sqrt\omega}{2\sqrt\omega}
	+\frac{-\epsilon_1(1-\epsilon_2)+d_2\sqrt\omega}{\sqrt\omega}
	\frac{\sqrt{\omega+(d_1/2)^2}-\sqrt\omega}{d_1},\\
	&\hat{L}^\s_{42}
	=\frac{-\epsilon_1\epsilon_2-d_2\sqrt s}{2\sqrt s}
	+\frac{-\epsilon_1\epsilon_2+d_2\sqrt s}{\sqrt s}\frac{\sqrt{s+(d_1/2)^2}-\sqrt s}{d_1},\\
	&\hat{L}^\u_{31}
	=\frac{\epsilon_1(1-\epsilon_2)-d_2\sqrt\omega}{2\sqrt\omega}
	+\frac{\epsilon_1(1-\epsilon_2)+d_2\sqrt\omega}{\sqrt\omega}\frac{\sqrt{\omega+(d_1/2)^2}-\sqrt\omega}{d_1},\\
	&\hat{L}^\u_{42}
	=\frac{\epsilon_1\epsilon_2-d_2\sqrt s}{2\sqrt s}
	+\frac{\epsilon_1\epsilon_2+d_2\sqrt s}{\sqrt s}\frac{\sqrt{s+(d_1/2)^2}-\sqrt s}{d_1}
\end{split}
\label{eqn:L}
\end{equation}
and
\[
\hat{L}_\pm=
\begin{pmatrix}
1&0&0&0&0&0\\
0&1&0&0&0&0\\
0&0&\pm\sqrt\omega&0&1&0\\
0&0&0&\pm\sqrt s&0&1
\end{pmatrix}.
\]
Here $\hat{L}_+$ and $\hat{L}_-$ are $4\times8$ real matrices consisting of
 bases in the subspaces spanned by row eigenvectors for $\D_{z,\eta}g(0,0;\beta_1,\epsilon_1,\epsilon_2)$
 such that the associated eigenvalues have negative and positive real parts, respectively.
In the computation, the distances $|z(x_\pm)|$ should be kept small
 but $|\eta(x_\pm)|$ do not have to be small necessarily
 like $|\zeta_{\R}(x_\pm)|,|\zeta_{\I}(x_\pm)|$ in Section~3.

Furthermore,
 to monitor the $L^2$ norm of $(\eta_1,\eta_2)^\T$,
 we add a parameter $c_1\in\Rset$ and the integral condition
\begin{align}
	\int_{x_-}^{x_+}\bigl(\eta_1(x)^2+\eta_2(x)^2\bigr)\d x=c_1.
	\label{eqn:ic3}
\end{align}
From \eqref{eqn:kernel} we see that
 if $(z(x),\eta(x))$ is a solution to \eqref{eqn:gker2} with $d_2=0$,
 then so is $(z(x),\eta(x)+\alpha\eta_0(t))$ for any $\alpha\in\Rset$,
 where $\eta_0(x)=(U'(x),V'(x),U''(x),V''(x))$.
To monitor the dependence of $\eta(x)$ on $\eta_0(x)$,
 we add a parameter $c_2\in\Rset$ and the integral condition
\begin{align}
	\int_{x_-}^{x_+}\bigl(\eta_1(x)z_3(x)+\eta_2(x)z_4(x)\bigr)\d x=c_2.
	\label{eqn:ic4}
\end{align}
If $c_2=0$, then $\eta(x)$ is orthogonal to $\eta_0(x)$.
To eliminate the multiplicity of solutions due to the translational symmetry of \eqref{eqn:SE}, we also add the integral condition \eqref{eqn:ic1}, as in Section~3.1.


\subsection{Numerical Results}

\begin{figure}[t]
	\includegraphics[scale=0.6]{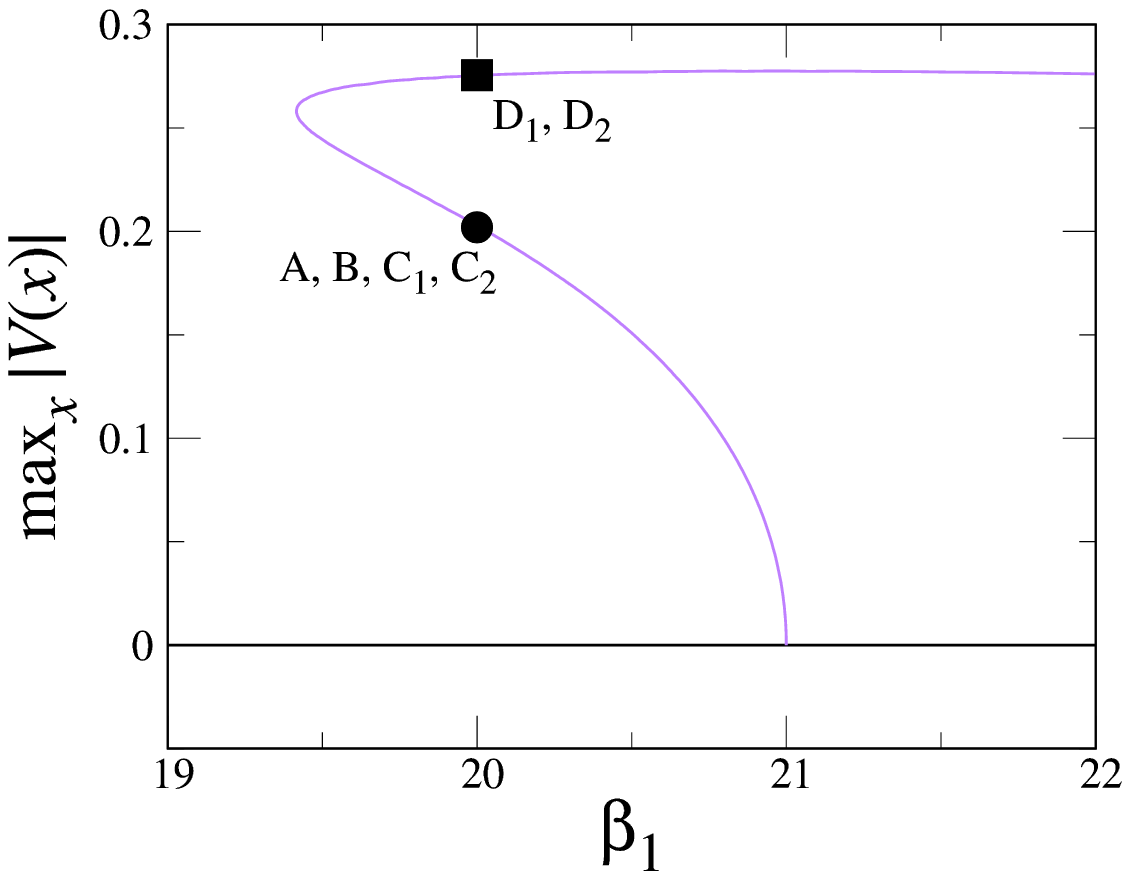}
	\caption{
		Partial enlargement of Fig.~\ref{fig:bif_diag}(a):
		Bifurcation diagram of solitary waves in \eqref{eqn:CNLS}.
		The labeled solutions in Table \ref{tab:runs} are located.
	}
	\label{fig:bif_diag_large}
\end{figure}

\begin{figure}[t]
	\begin{minipage}{0.495\textwidth}
		\includegraphics[scale=0.5]{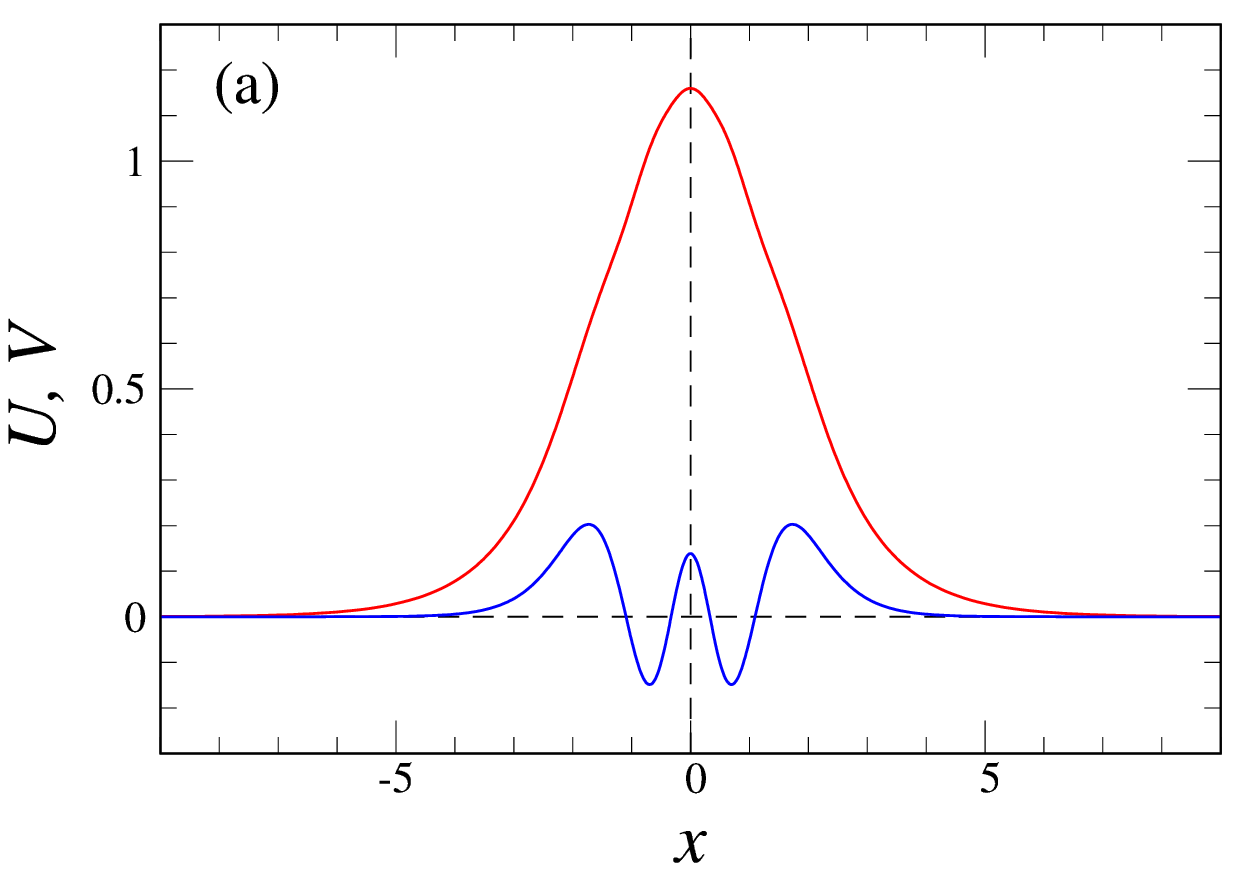}
	\end{minipage}
	\begin{minipage}{0.495\textwidth}
	\begin{center}
		\includegraphics[scale=0.5]{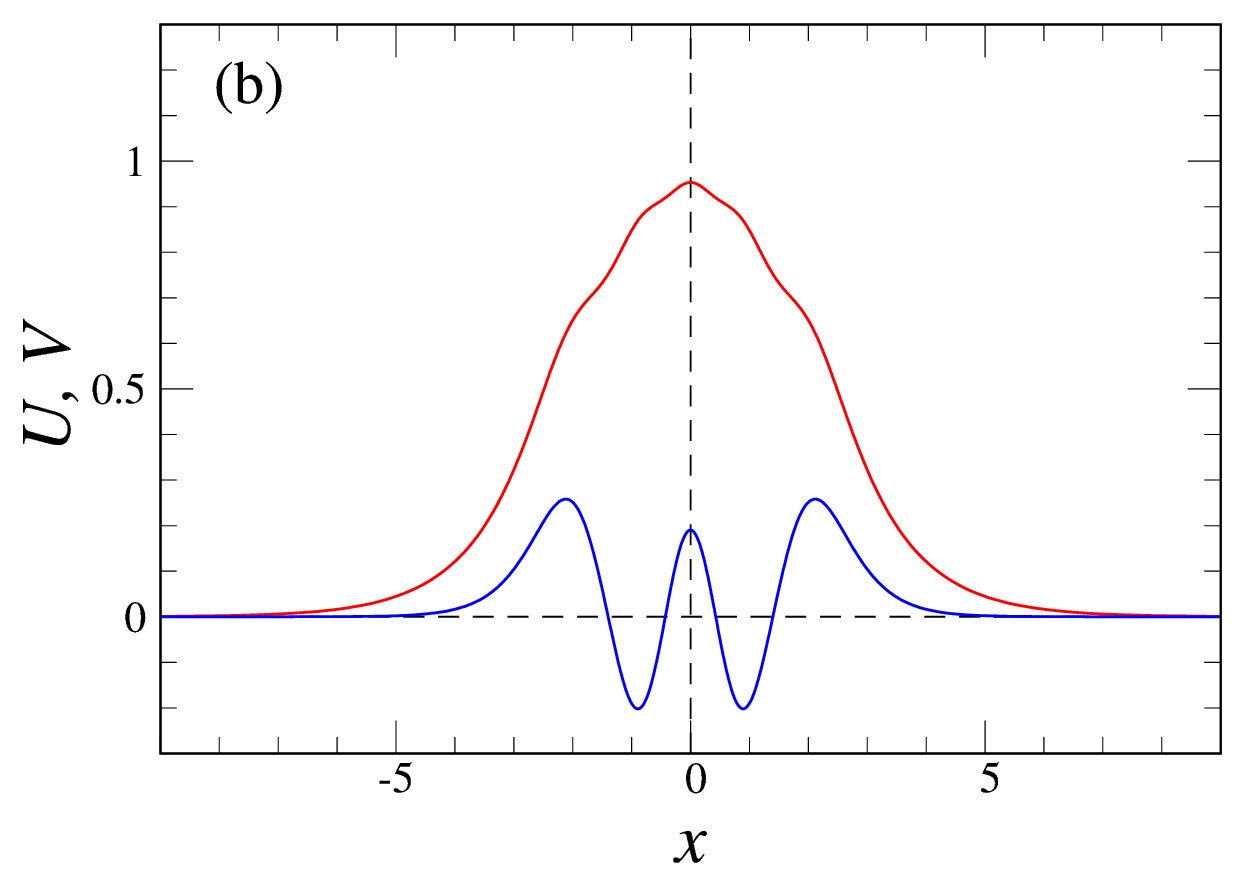}
	\end{center}
	\end{minipage}
	\begin{center}
		\includegraphics[scale=0.5]{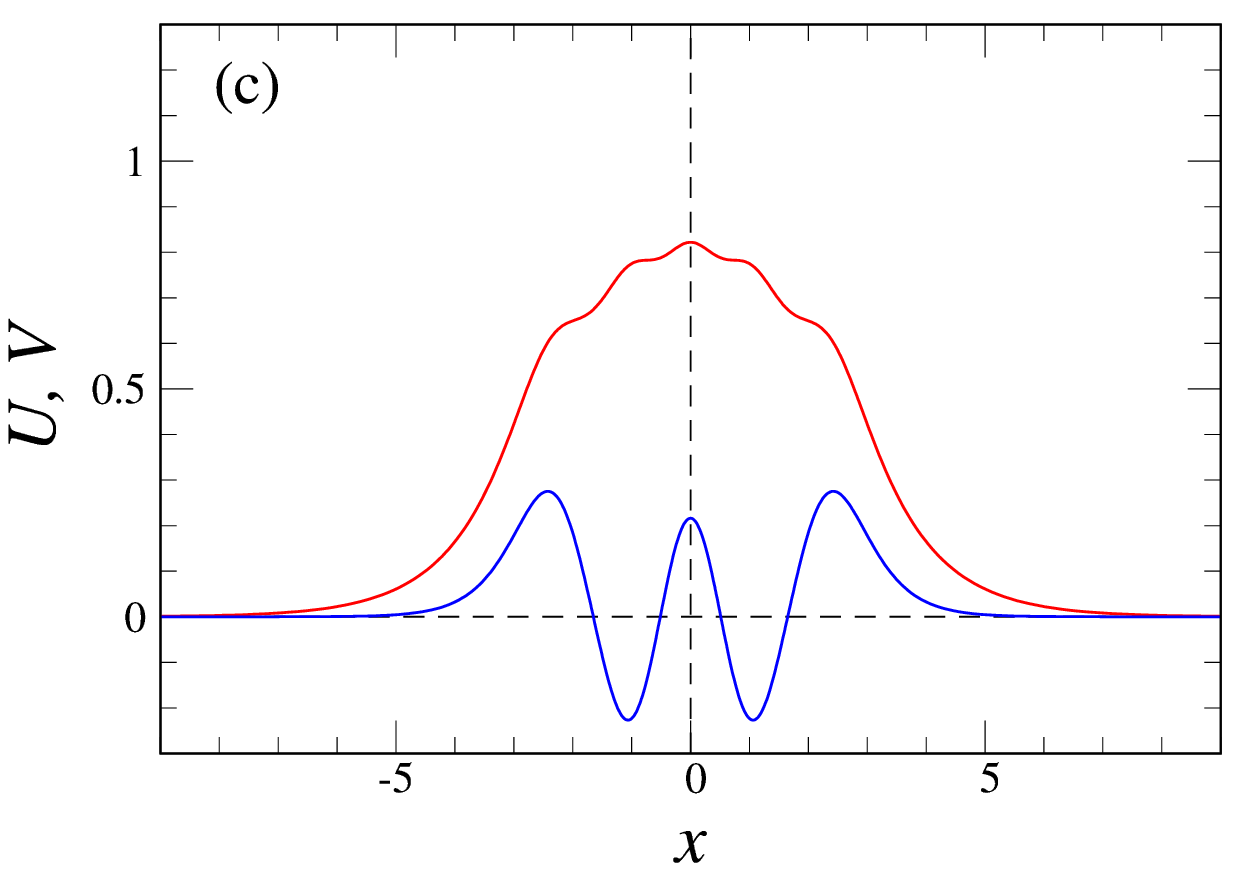}
	\end{center}
	\caption{
		Numerically computed homoclinic solutions to \eqref{eqn:SE}
		for $(\omega,s,\beta_2)=(1,4,2)$:
		(a) $\beta_1=20$;
		(b) $\beta_1=\beta_1^\mathrm{SN}$; 
		(c) $\beta_1=20$.
		The red and blue lines, respectively,
		represents the $U$- and $V$-components.
		The solutions in plates~(a) and (c) are, respectively,
		 denoted by $\bullet$ and {\tiny$\blacksquare$}
		 on the branch in Fig.~\ref{fig:bif_diag_large}.
	}
	\label{fig:gker_UV}
\end{figure}

We set $\omega=1$, $s=4$, and $\beta_2=2$ as in Section 3.2.
Figure \ref{fig:bif_diag_large} is a partial enlargement of Fig.~\ref{fig:bif_diag}(a)
 in which the pitchfork and saddle-node bifurcation points on the fifth branch
 are contained.
The corresponding homoclinic solutions to \eqref{eqn:SE} on the branch
 are displayed in Fig.~\ref{fig:gker_UV}.
We carried out numerical computations stated in Section~4.1
 along the branch beyond the saddle-node bifurcation point
 $\beta_1=\beta_1^\mathrm{SN}\approx19.41626$.
We also used the computer tool \texttt{AUTO} \cite{D12}
 to obtain numerical solutions to \eqref{eqn:gker2}
 satisfying the boundary condition \eqref{eqn:bc3}
 under the integral conditions \eqref{eqn:ic1}, \eqref{eqn:ic3}, and \eqref{eqn:ic4},
 as in Section~3.

\begin{table}[t]
	\caption{Summary of the numerical continuations.
	Here $c_0\approx0.074836$ (see \eqref{eqn:c0}).}\label{tab:runs}
\renewcommand{\arraystretch}{1.4}
	\begin{tabular}{|c|c|c|c|c|c|c|c|} \hline
		Run&Varied parameter&\multicolumn{4}{|c|}{Fixed parameter values}& Starting&Terminating\\ \cline{3-6}
		no.&&\makebox[1.8em]{$\beta_1$}&\makebox[1.8em]{$\epsilon_2$}&\makebox[1.8em]{$c_1$}&\makebox[1.8em]{$c_2$} &solution&solution\\ \hline\hline
		1&$c_1=c_0\to1$&$20$&$0$&-&$c_0$&A&B\\ \hline
		2&$c_2=c_0\to0$&$20$&$0$&$1$&-&B&$\mathrm{C}_1$\\ \hline
		3&$\beta_1=20\to\beta_1^\mathrm{SN}\to20$&-&$0$&$1$&$0$&$\mathrm{C}_1$&$\mathrm{D}_1$\\ \hline
		4&$\epsilon_2=0\to1$&$20$&-&$1$&$0$&$\mathrm{C}_1$&$\mathrm{C}_2$\\ \hline
		5&$\beta_1=20\to\beta_1^\mathrm{SN}\to20$&-&$1$&$1$&$0$&$\mathrm{C}_2$&$\mathrm{D}_2$\\ \hline
	\end{tabular}
\end{table}

Since at the fifth pitchfork bifurcation point $\beta_1=\beta_1^{(4)}=21$,
 $\dim\Ker J\mathcal L$ increases by two
 and consequently \eqref{eqn:gker2} is highly degenerate,
 it is difficult to continue a branch of solutions
 in the boundary value problem beyond there.
From this reason
 we computed the solution branch after the pitchfork bifurcation occurs.
As the starting solution in a series of numerical continuations, we adopted
\begin{align}
	z=(U(x),V(x),U'(x),V'(x))^\T\quad
	\eta=(U'(x),V'(x),U''(x),V''(x))^\T
	\label{eqn:start_sol}
\end{align}
at $\beta_1=20$ on the branch
 with $\epsilon_1,\epsilon_2,d_1,d_2=0$,
 where $U(x),V(x),U'(x)$, and $V'(x)$ were numerically obtained along with
\begin{align*}
&
U''(x)=\omega U(x)-(U(x)^2+\beta_1V(x)^2)U(x),\\
&
V''(x)=s V(x)-(\beta_1U(x)^2+\beta_2V(x)^2)V(x)
\end{align*}
in advance by another numerical continuation for the boundary value problem
 of \eqref{eqn:dyn_sys} with \eqref{eqn:bc1}.
We chose $x_\pm=\pm9$ and executed five runs in total.
All of the runs are summarized in Table \ref{tab:runs}.
In the numerical continuations,
 $\beta_1$, $\epsilon_2$, $c_1$ or $c_2$ was varied
 while $\epsilon_1$, $d_1$, and $d_2$ were taken as the free parameters.
The distances $|z(x_\pm)|$ 
 were monitored along with $|\eta(x_\pm)|$
 and kept small ($\approx10^{-3}$) 
 during the computations.
Moreover, since $d_1$ is very small (it should be zero theoretically), the approximations
\[
\frac{\sqrt{\omega+(d_1/2)^2}-\sqrt{\omega}}{d_1}
\approx\frac{d_1}{8\sqrt{\omega}},\quad
\frac{\sqrt{s+(d_1/2)^2}-\sqrt{s}}{d_1}
\approx\frac{d_1}{8\sqrt{s}},
\]
were used in the computations of \eqref{eqn:L}.

In the first run, we took the numerical solution \eqref{eqn:start_sol}
 with $(\beta_1,\epsilon_2)=(20,0)$ and
\begin{align}
	c_1=c_2=\int_{x_-}^{x_+}\bigl(U'(x)^2+V'(x)^2\bigr)\,\d x
	=:c_0\approx0.074836
\label{eqn:c0}
\end{align}
labeled by `A' as the starting solution
 and continued it from $c_1=c_0$ to $1$ for $\beta_1,\epsilon_2,c_2$ fixed.
The solution calculated at $c_1=1$,
 the $(\eta_1,\eta_2)$-components of which are normalized,
 is labeled by `B'.
In the second run, we fixed $(\beta_1,\epsilon_2,c_1)=(20,0,1)$
 and followed the solution `B' from $c_2=c_0$ to $0$.
The solution calculated at $c_2=0$,
 for which $(\eta_1(x),\eta_2(x))$ is orthogonal to $(z_3,z_4)=(U'(x),V'(x))$ by \eqref{eqn:ic4},
 is labeled by `$\mathrm{C}_1$'.
In the third run,
 we fixed $(\epsilon_2,c_1,c_2)=(0,1,0)$ and followed the solution `$\mathrm{C}_1$' for $J\mathcal L\psi=\epsilon_1\varphi_2$ (see \eqref{eqn:gker})
 from $\beta_1=20$ to $\beta_1^\mathrm{SN}$
 and from $\beta_1^\mathrm{SN}$ to $20$.
The finally obtained solution is labeled by `$\mathrm{D}_1$'.
See Fig.~\ref{fig:bif_diag_large}.

\begin{figure}[t]
	\begin{minipage}{0.495\textwidth}
		\includegraphics[scale=0.5]{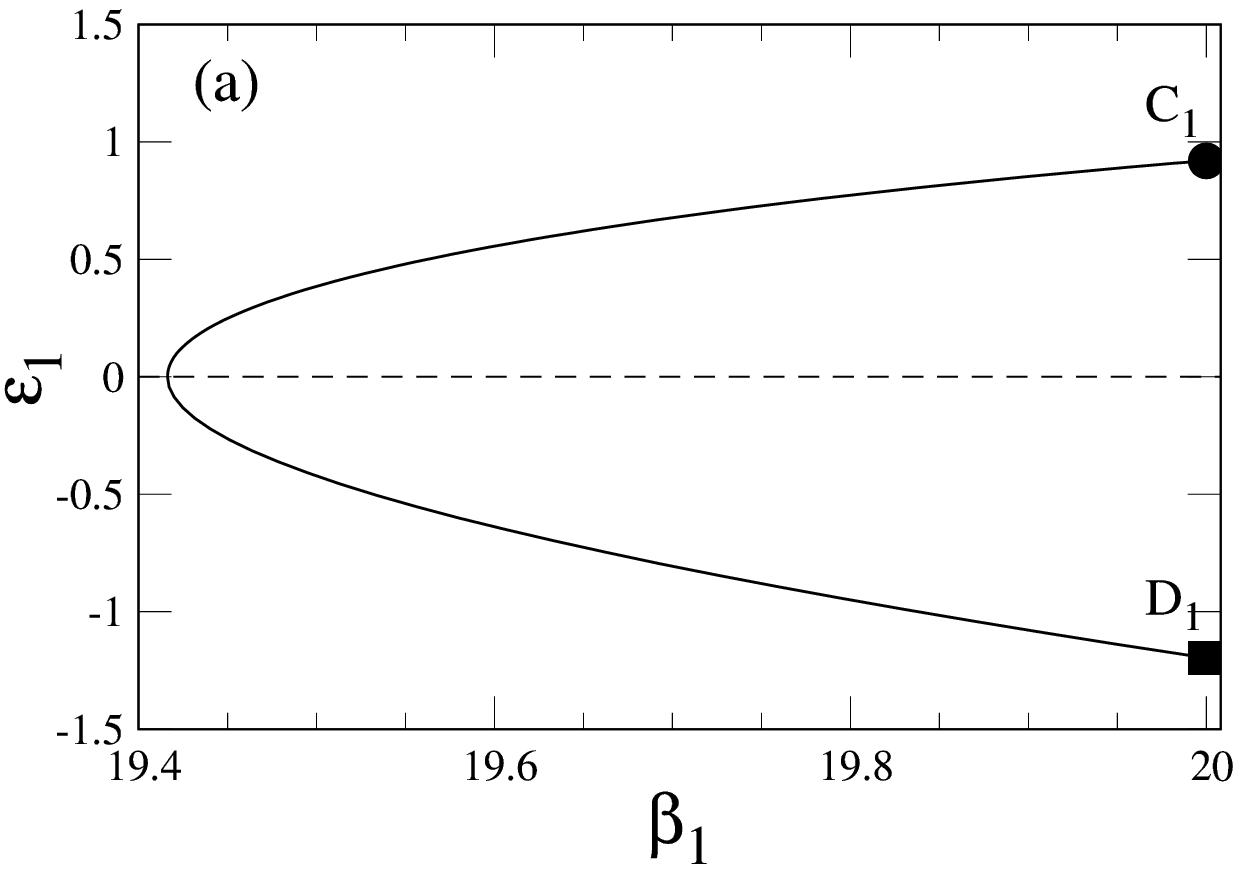}
		\includegraphics[scale=0.5]{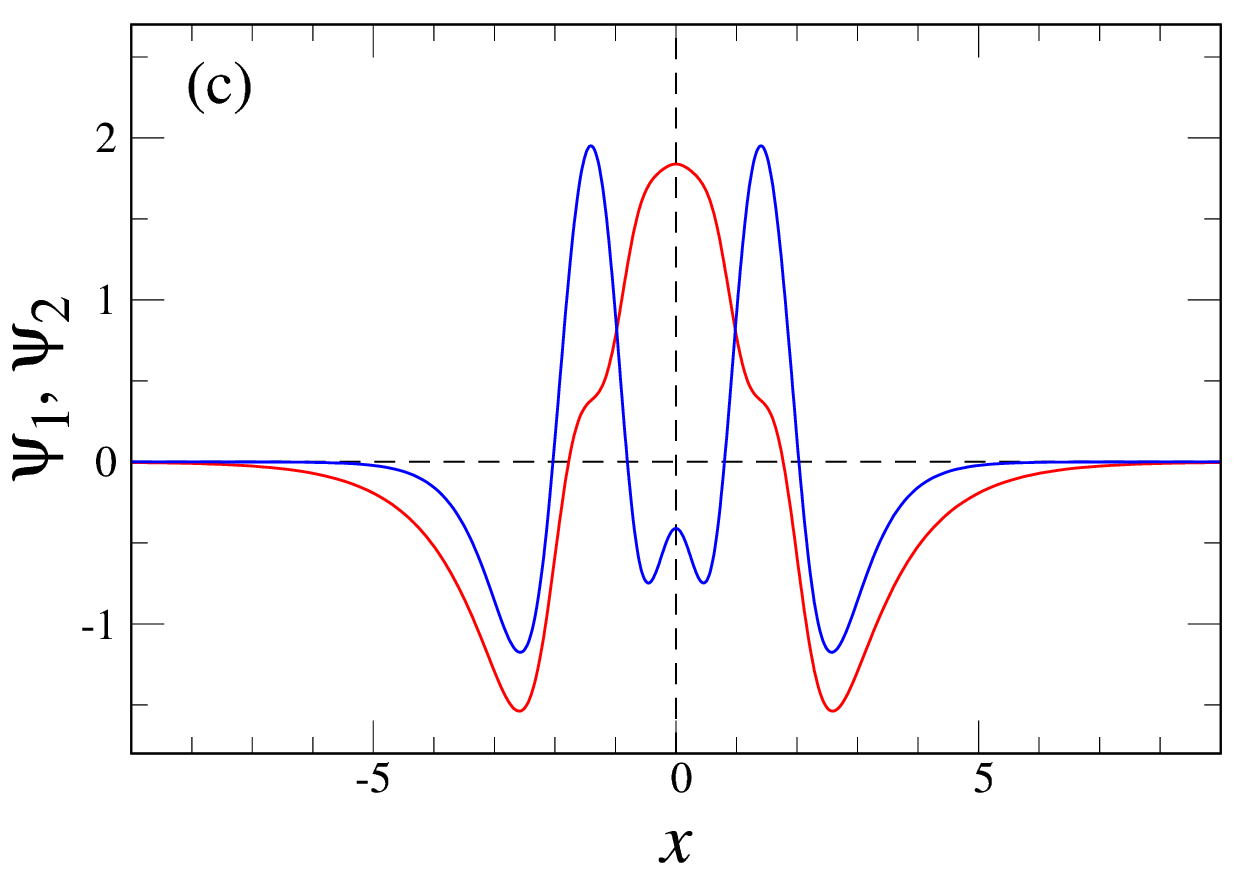}
	\end{minipage}
	\begin{minipage}{0.495\textwidth}
	\begin{center}
		\includegraphics[scale=0.5]{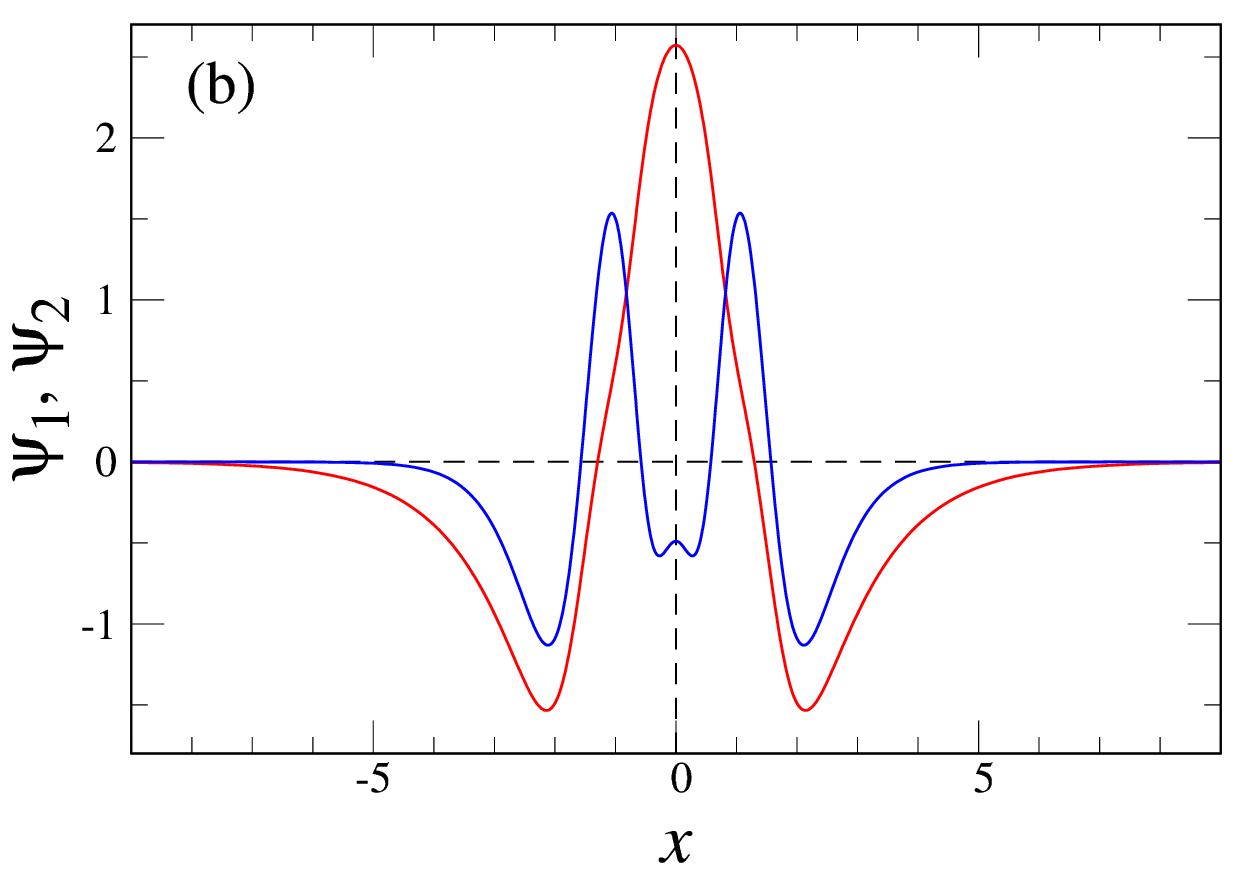}
		\includegraphics[scale=0.5]{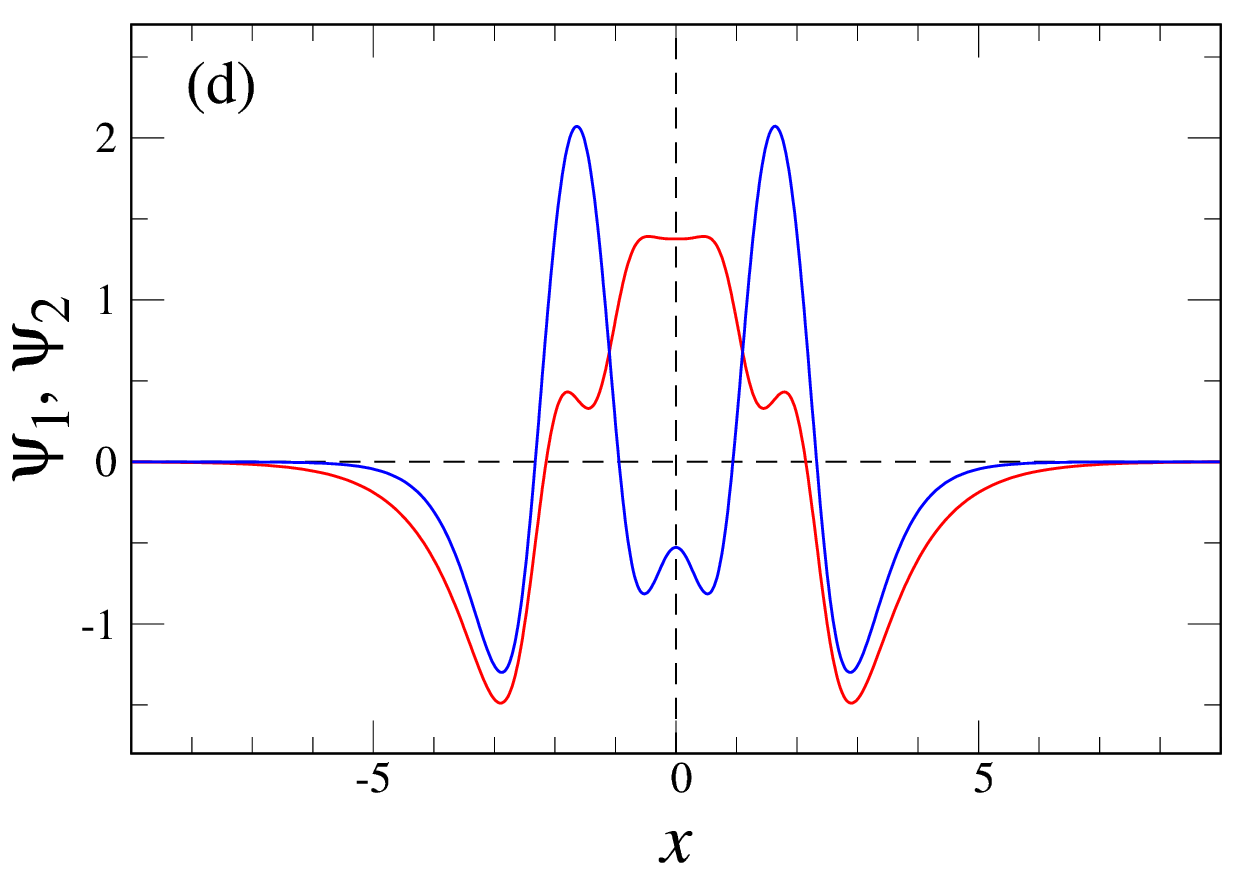}
	\end{center}
	\end{minipage}
	\caption{
		Numerically computed generalized eigenfunctions satisfying \eqref{eqn:gker}
		with $\epsilon_2=0$ for $(\omega,s,\beta_2)=(1,4,2)$:
		(a) Solution branch $\mathrm{C}_1\mathrm{D}_1$;
		(b) $\beta_1=20$ (at $\mathrm{C}_1$);
		(c) $\beta_1=\beta_1^\mathrm{SN}$;
		(d) $\beta_1=20$ (at $\mathrm{D}_1$).
		The red and blue lines, respectively, represent
		the first and second components.
	}
	\label{fig:gker0}
\end{figure}

Figure~\ref{fig:gker0} shows generalized eigenfunctions satisfying \eqref{eqn:gker}
 with $\epsilon_2=0$ along with the solution branch obtained from the third run.
In particular, at the saddle-node bifurcation point $\beta_1=\beta_1^\mathrm{SN}$,
 we observe $\epsilon_1=0$, so that the generalized eigenfunction 
 expressed as a linear combination of $\varphi_1$ and $\chi_2$
 becomes an eigenfunction for the zero eigenvalue
 in the eigenvalue problem \eqref{eqn:eig_eq}, as we suspect.
Thus, the corresponding solitary wave does not change its stability type
 at $\beta_1=\beta_1^\mathrm{SN}$: the number of eigenvalues with positive real parts does not change.
This is similar to Yang's result \cite{Ya12b}
 for general single NLS equations with external potentials.

In the fourth run, we fixed $(\beta_1,c_1,c_2)=(20,1,0)$
 and followed the solution `$\mathrm{C}_1$' from $\epsilon_2=0$ to $1$
 for $\beta_1=20$ fixed.
The solution calculated at $\epsilon_2=1$ is labeled by `$\mathrm{C}_2$'.
In the last run, we fixed $(\epsilon_2,c_1,c_2)=(1,1,0)$
 and followed the solution `$\mathrm{C}_2$' {for $J\mathcal L\psi=\epsilon_1\varphi_3$ (see \eqref{eqn:gker})}
 from $\beta_1=20$ to $\beta_1^\mathrm{SN}$
 and from $\beta_1^\mathrm{SN}$ to $20$.
The finally obtained solution is labeled by `$\mathrm{D}_2$'.
See Fig.~\ref{fig:bif_diag_large}.

\begin{figure}[t]
	\begin{minipage}{0.495\textwidth}
		\includegraphics[scale=0.5]{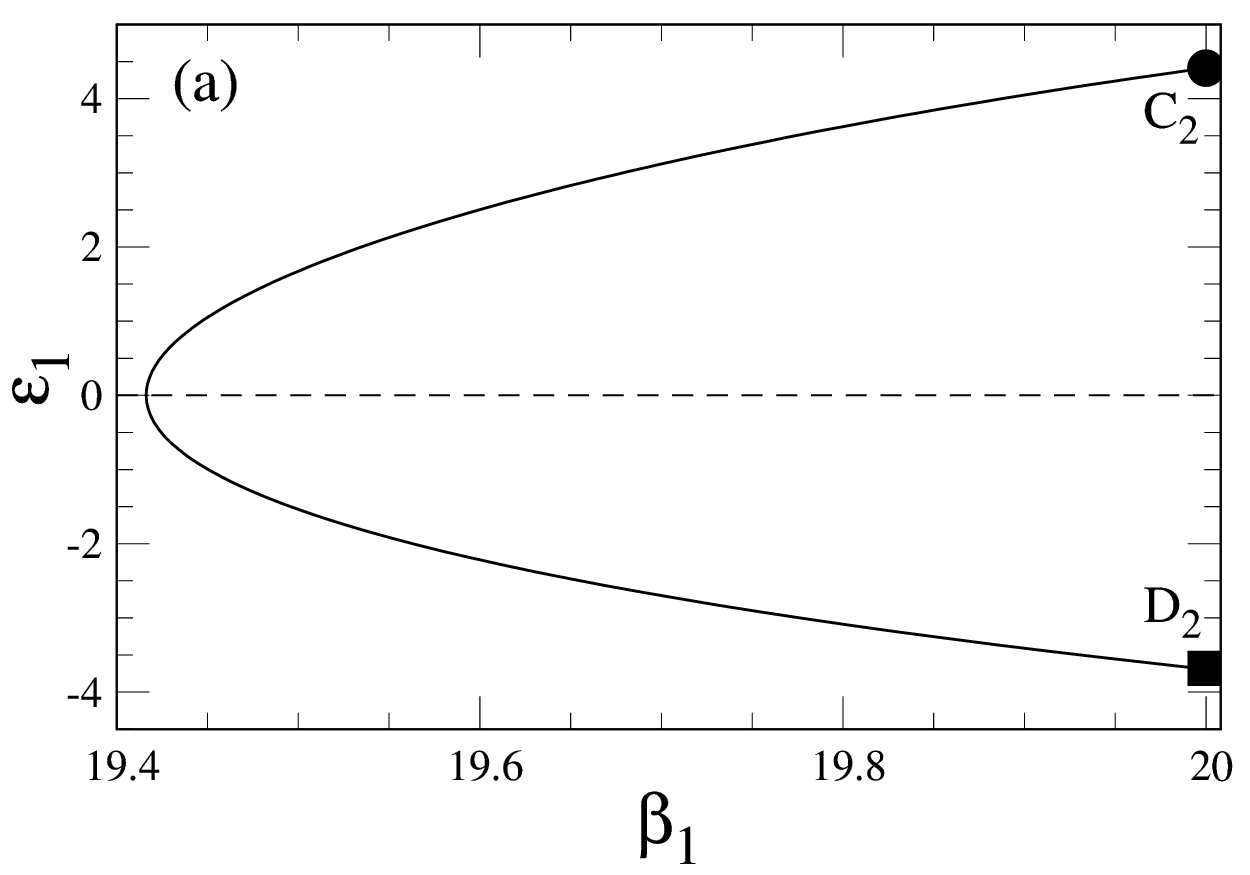}
		\includegraphics[scale=0.5]{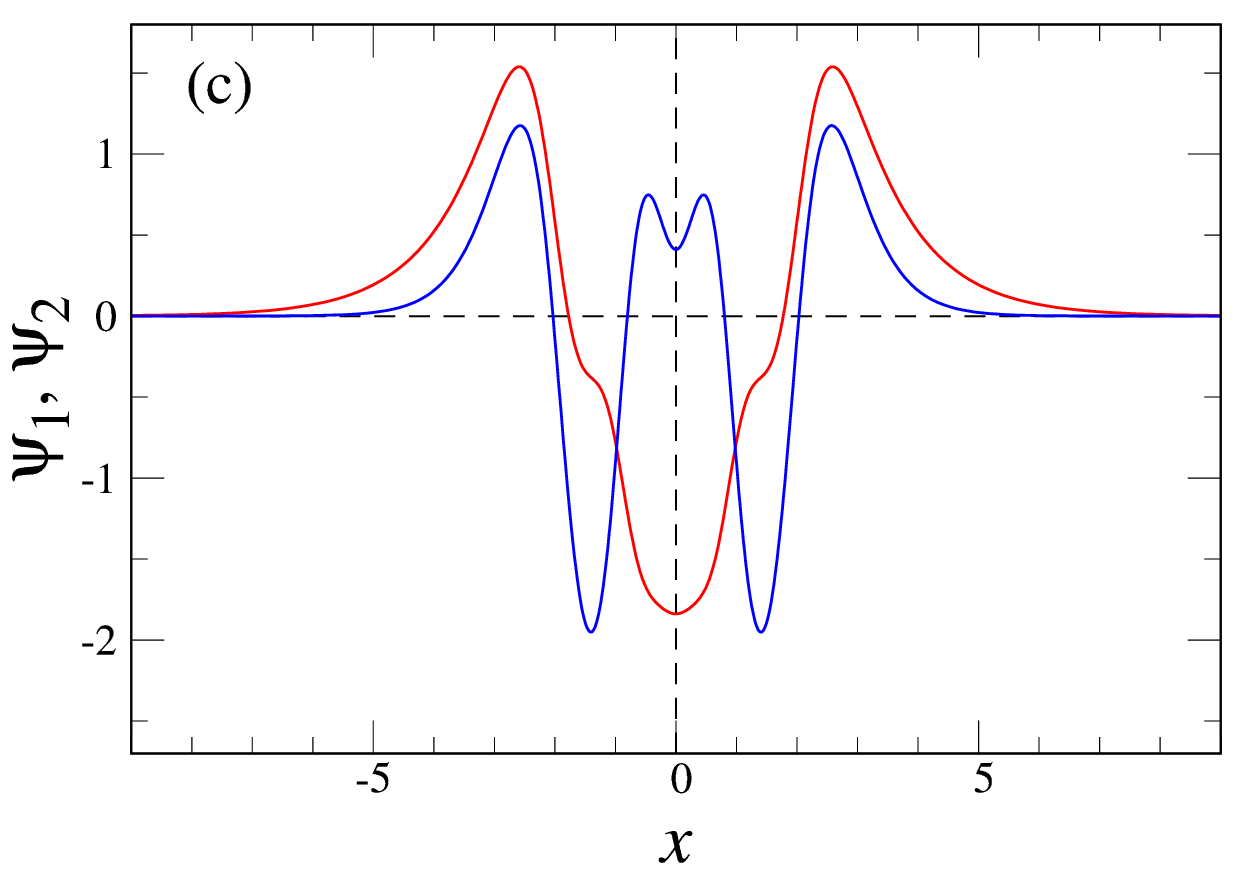}
	\end{minipage}
	\begin{minipage}{0.495\textwidth}
	\begin{center}
		\includegraphics[scale=0.5]{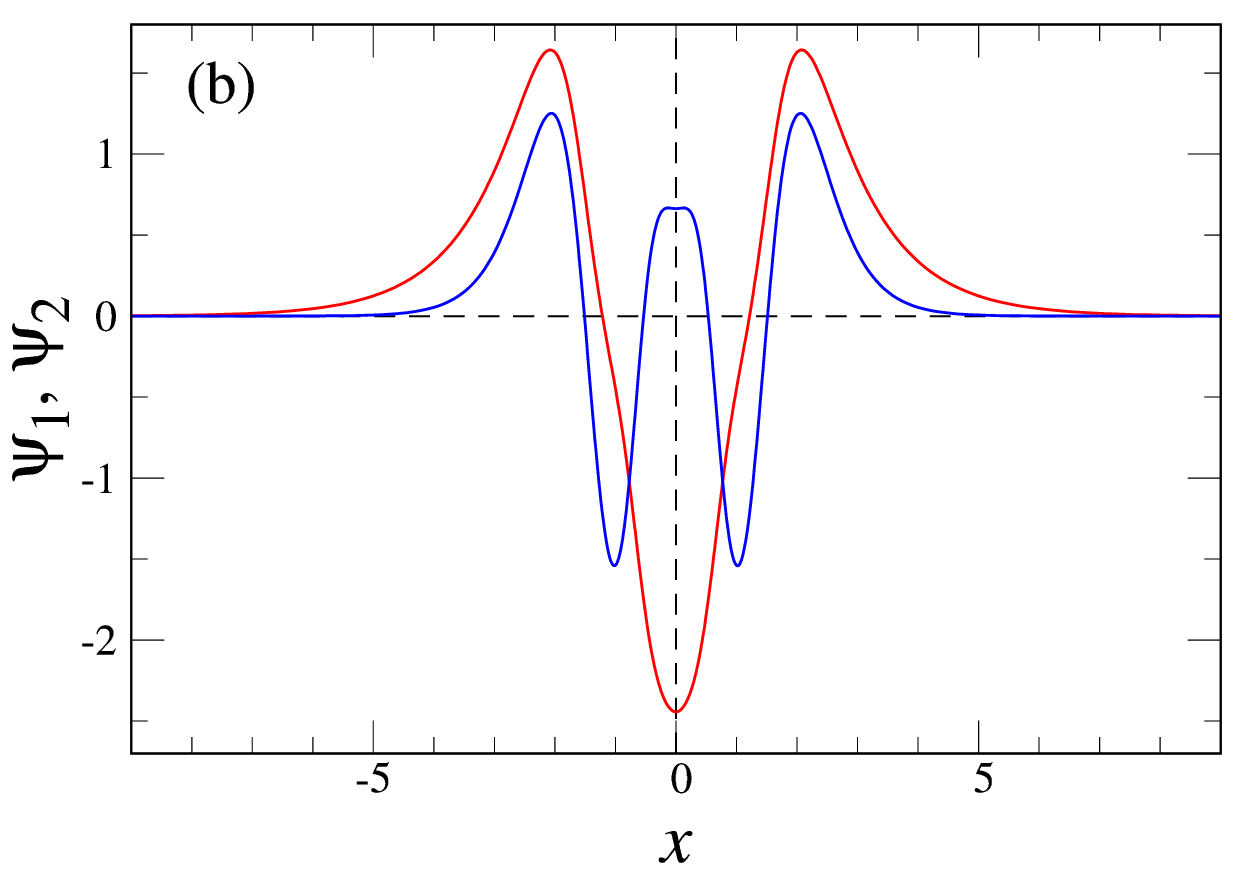}
		\includegraphics[scale=0.5]{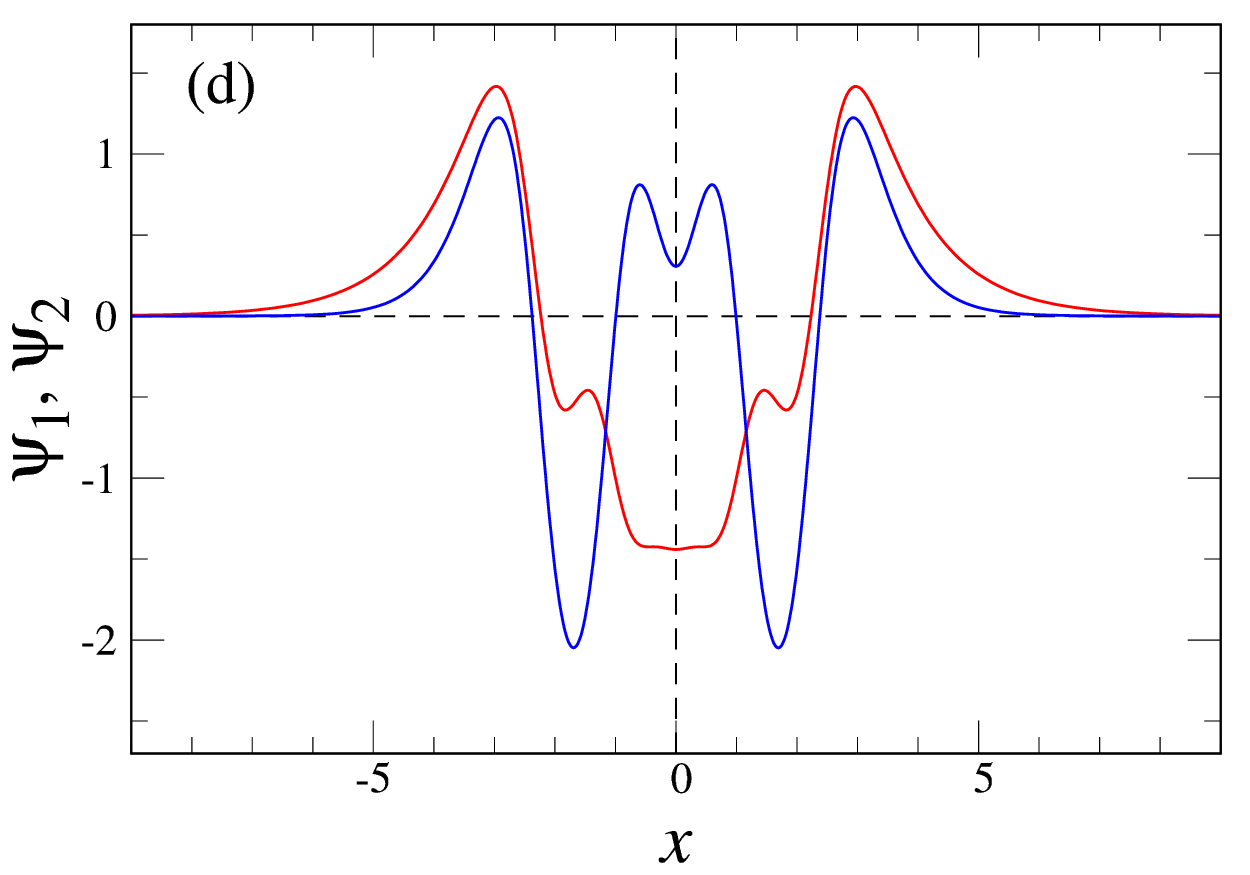}
	\end{center}
	\end{minipage}
	\caption{
		Numerically computed generalized eigenfunctions satisfying \eqref{eqn:gker}
		with $\epsilon_2=1$ for $(\omega,s,\beta_2)=(1,4,2)$:
		(a) Solution branch $\mathrm{C}_2\mathrm{D}_2$;
		(b) $\beta_1=20$ (at $\mathrm{C}_2$);
		(c) $\beta_1=\beta_1^\mathrm{SN}$;
		(d) $\beta_1=20$ (at $\mathrm{D}_2$).
		The red and blue lines, respectively, represent
		the first and second components.
	}
	\label{fig:gker1}
\end{figure}

Figure~\ref{fig:gker1} shows generalized eigenfunctions satisfying \eqref{eqn:gker}
 with $\epsilon_2=1$ along with the solution branch obtained from the last run.
In particular, at the saddle-node bifurcation point $\beta_1=\beta_1^\mathrm{SN}$,
 we observe $\epsilon_1=0$, so that the generalized eigenfunction 
 expressed as a linear combination of $\varphi_1$ and $\chi_3$
 becomes an eigenfunction for the zero eigenvalue
 in the eigenvalue problem \eqref{eqn:eig_eq}, again.
Moreover, the eigenfunction of Fig.~\ref{fig:gker1}(c)
 coincide with that of Fig.~\ref{fig:gker0}(c) up to multiplication by $-1$.
Thus, the generalized eigenfunctions for $\varphi_2$ and $\varphi_3$ give the same eigenfunction at $\beta_1=\beta_1^{\mathrm{SN}}$, and $\dim\Ker J\mathcal L$ increases by one at $\beta_1=\beta_1^\mathrm{SN}$.

We close this paper with showing that $\dim\gKer J\mathcal L$ does not change
 at $\beta_1^\mathrm{SN}$
 although the two linearly independent generalized eigenfunctions were observed 
 to converge to the eigenfunction as $\beta_1\to\beta_1^\mathrm{SN}$.

Fix $\beta_1=\beta_1^\mathrm{SN}$.
Let $(\eta_1^\mathrm{SN},\eta_2^\mathrm{SN})$ denote
 the $(\eta_1,\eta_2)$-components of the solution to \eqref{eqn:gker2}
 with $(\epsilon_1,c_2)=(0,0)$ 
 corresponding to an eigenfunction of \eqref{eqn:eig_eq} for the zero eigenvalue,
 such as plotted in Figs.~\ref{fig:gker0}(c) and \ref{fig:gker1}(c).
We see that $(\eta_1^\mathrm{SN},\eta_2^\mathrm{SN})\in\Ker\mathcal L_+$,
 where $\mathcal L_+$ is the linear operator given in \eqref{eqn:def_Lpm}.
Moreover, $\Ker\mathcal L_+$ is of dimension two at most
 since the corresponding four-dimensional system of first-order ODEs
 converges to \eqref{eqn:dyn_sys0} as $x\to\pm\infty$
 and the stable and unstable subspaces of the origin in \eqref{eqn:dyn_sys0}
 are of dimension two
 ($\dim E^\s=\dim E^\u=2$).
So $\{(U',V'),(\eta_1^\mathrm{SN},\eta_2^\mathrm{SN})\}$
 is a basis of $\Ker\mathcal L_+$.
Similarly, $\Ker\mathcal L_-$ is of dimension two at most
 and $\{(U,0),(0,V)\}$ is its basis,
 where $\mathcal L_-$ is the linear operator given in \eqref{eqn:def_Lpm}.
Thus, $\Ker J\mathcal L=\Ker\mathcal L_+\oplus\Ker\mathcal L_-$ is of dimension four
 and spanned by
\begin{align*}
	\{\varphi_1,(\eta_1^\mathrm{SN},\eta_2^\mathrm{SN},0,0)^\T,
	\varphi_2, 
	\varphi_3\},
\end{align*}
where $\varphi_j$, $j=1,2,3$, were given in \eqref{eqn:kernel}.

To determine $\gKer J\mathcal L$, we consider the solvability of
\begin{align}
	\mathcal L_-
	\begin{pmatrix}
		\zeta_1\\
		\zeta_2
	\end{pmatrix}
	=\alpha_1\begin{pmatrix}
		U'\\
		V'
	\end{pmatrix}
	+\alpha_2\begin{pmatrix}
		\eta_1^\mathrm{SN}\\
		\eta_2^\mathrm{SN}
	\end{pmatrix}
	\label{eqn:gker_eq1}
\end{align}
and
\begin{align}
	\mathcal L_+
	\begin{pmatrix}
		\zeta_1\\
		\zeta_2
	\end{pmatrix}
	=\alpha_1\begin{pmatrix}
		U\\
		0
	\end{pmatrix}
	+\alpha_2\begin{pmatrix}
		0\\
		V
	\end{pmatrix},
	\label{eqn:gker_eq2}
\end{align}
where $\alpha_j$, $j=1,2$, are constants.
Note that nontrivial solutions to \eqref{eqn:gker_eq1} and \eqref{eqn:gker_eq2}
 provide elements of $\gKer J\mathcal L$.
Numerical integrations carried out in the software \texttt{AUTO} yielded
\begin{align*}
	I_1=\int_{x_-}^{x_+}\eta_1^\mathrm{SN}(x)U(x)\,\d x\approx1.492,\quad
	I_2=\int_{x_-}^{x_+}\eta_2^\mathrm{SN}(x)V(x)\,\d x\approx-0.373,
	\label{eqn:inn_val}
\end{align*}
which indicates along with  the Fredholm alternative theorem \cite{KaP13} that
 there exists an $L^2$ solution to \eqref{eqn:gker_eq1} (resp.\ to \eqref{eqn:gker_eq2})
 if and only if $\alpha_2=0$ (resp.\ $I_1\alpha_1+I_2\alpha_2=0$).
Actually, $\zeta=(-xU(x)/2,-xV(x)/2)$ is the solution
 to \eqref{eqn:gker_eq1} with $\alpha_2=0$,
 and consequently it is reconfirmed that $\chi_1\in\gKer J\mathcal L$.
Thus, there exist two linearly independent generalized eigenfunctions 
 and $\dim\gKer J\mathcal L$ does not change at  $\beta_1=\beta_1^\mathrm{SN}$.
 

\end{document}